\documentclass[a4paper,12pt]{article}

\usepackage{latexsym,makeidx}
\usepackage{mathrsfs}
\usepackage{amsfonts,amsmath,amssymb}
\usepackage{indentfirst}
\usepackage{subeqnarray}  
\usepackage{verbatim}
\usepackage[pdftex]{graphicx}
\usepackage{epstopdf}

\usepackage[mathscr]{eucal}
\usepackage{graphicx, color}

\newtheorem{Theorem}{Theorem}

\newtheorem{Lemma}[Theorem]{\sc Lemma}

\newtheorem{Remark}[Theorem]{\sc Remark}

\newtheorem{Problem}[Theorem]{\sc Problem}

\newcommand{\R}{{\if mm {\rm I}\mkern -3mu{\rm R}\else \leavevmode
\hbox{I}\kern -.17em\hbox{R} \fi}}

\newcommand{\bu}{\mbox{\boldmath{$u$}}}
\newcommand{\bv}{\mbox{\boldmath{$v$}}}
\newcommand{\bw}{\mbox{\boldmath{$w$}}}
\newcommand{\bx}{\mbox{\boldmath{$x$}}}

\newcommand{\fb}{\mbox{\boldmath{$f$}}}

\newcommand{\bsigma}{\mbox{\boldmath{$\sigma$}}}

\newcommand{\bvarepsilon}{\mbox{\boldmath{$\varepsilon$}}}
\newcommand{\bnu}{\mbox{\boldmath{$\nu$}}}

\newcommand{\bzero}{\mbox{\boldmath{$0$}}}

\newcommand{\bz}{\mbox{\boldmath{$z$}}}

\newcounter{theorem}

\def\sqr#1#2{{
    \vcenter{
         \vbox{\hrule height.#2pt
               \hbox{\vrule width.#2pt height#1pt \kern#1pt
                     \vrule width.#2pt
               }
               \hrule height.#2pt
         }
    }
}}

\def\bar{\overline}
\def\real{\mathbb{R}}

\def\lista#1
{{ \itemindent 0.0cm \labelsep .2cm \leftmargin 0.8cm \rightmargin
0.0cm \labelwidth 0.6cm \topsep 0.0mm
\parsep 0.0mm
\itemsep 0.0mm
\begin{list}{}
{ \setlength{\leftmargin}{.8cm} \setlength{\rightmargin}{0.0cm}
\setlength{\parsep}{0.0mm} \setlength{\topsep}{.0mm}
\setlength{\parskip}{.0cm} \setlength{\itemsep}{.0cm} }
{#1}\end{list}} }

 \numberwithin{equation}{section}
\begin{document}
	
\title{
A new class of history-dependent quasi variational-hemivariational inequalities with constraints
\thanks{\,
Project was supported by the European Union's Horizon 2020 Research and Innovation Programme under the Marie Sk{\l}odowska-Curie grant agreement No. 823731 CONMECH, the NNSF of China Grant Nos. 12001478 and 12101143, the Natural Science Foundation of Guangxi Grants Nos. 2021GXNSFFA196004, 2020GXNSFBA297137, GKAD21220144, 2018GXNSFAA281353 and 2019GXNSFBA185005),
the Beibu Gulf University Project No. 2018KYQD06, the Startup Project of Doctor Scientific Research of Yulin Normal University No. G2020ZK07,
the Mini\-stry of Science and Higher Education of Republic of Poland under Grants Nos. 4004/GGPJII/H2020/2018/0 and 440328/PnH2/2019,
and the National Science Centre of Poland under Project No. 2021/41/B/ST1/01636.
}}
	
\author{
Stanis{\l}aw Mig\'orski $^{1,3}$
\footnote{\,
E-mail address: stanislaw.migorski@uj.edu.pl.}, \ \
Yunru Bai $^2$
\footnote{\,
E-mail address: yunrubai@163.com.
}
\ \ and \ \
Shengda Zeng  $^{4}$
\footnote{\, E-mail address: zengshengda@163.com; shengdazeng@gmail.com.  Corresponding author.}
}
	
\date{}
\maketitle

\begin{center}
$^1$ College of Applied Mathematics\\
 Chengdu University of Information Technology\\
  Chengdu, 610225, Sichuan Province, P.R. China\\[2mm]
 $^2$ School of Science \\
Guangxi University of Science and Technology\\
 Liuzhou 545006, Guangxi Province, P.R. China\\[2mm]
$^3$ Jagiellonian University in Krakow\\
Faculty of Mathematics and Computer Science\\
ul. Lojasiewicza 6, 30348 Krakow, Poland\\[2mm]
$^4$  Guangxi Colleges and Universities Key Laboratory\\
 of Complex System Optimization and Big Data Processing\\
   Yulin Normal University\\
    Yulin 537000, Guangxi Province, P.R. China
\end{center}

\thispagestyle{empty}


\begin{abstract}
In this paper we consider an abstract class
of time-dependent quasi variational-hemivariational inequalities which involves history-dependent operators and a set of unilateral constraints.
First, we establish the existence and uniqueness
of solution by using a recent result for elliptic variational-hemivariational inequalities in reflexive Banach spaces combined with a fixed-point principle for history-dependent operators.
Then, we apply the abstract result
to show the unique weak solvability to a quasistatic   viscoelastic frictional contact problem.
The contact law involves a unilateral Signorini-type condition for the normal velocity
and the nonmonotone normal damped response condition while the friction condition is a version of the Coulomb law of dry friction in which the friction bound depends on the accumulated slip.
\end{abstract}

\smallskip


\smallskip










\section{Introduction}\label{Introduction}

\noindent
The theory of hemivariational and variational--hemivariational inequalities has been ori\-ginated in early 1980s with the pio\-neering works of Panagiotopoulos, see~\cite{NP,P2,P}.
These inequalities have served as variational descriptions of many physical phenomena that include nonconvex, nondifferentiable and locally Lipschitz potentials, and they have played an important role in a description of diverse mechanical problems arising in solid and fluid mechanics.
The theory has undergone a remarkable development in pure
and applied mathematics, see monographs~\cite{CLM1,KS,MOSBook,NP,P} for the
mathematical theory, applications, and related issues.

The aim of this paper is to study a general class
of time-dependent quasi variational-hemivariational inequalities involving history-dependent operators and a set of constraints.
Given the operators $A$, $f$ and $M$,
a convex potential $\varphi$,
a locally Lipschitz potential $j$,
history-dependent operators $R_1$, $R_2$, $R_3$
and $R_4$, and a constraint set $K$, the problem under investigation reads as follows:
find $w \in L^2(0, T; K)$ such that
for a.e. $t \in (0, T)$, we have
\begin{equation*}
\begin{cases}
\langle A (t, (R_1 w)(t), w(t)) - f(t, (R_2 w)(t)),
v - w(t) \rangle \\[1mm]
\qquad + \,
\varphi (t, (R_3 w)(t), w(t), v)
- \varphi(t, (R_3 w)(t), w(t), w(t)) \\[1mm]
\ \ \ \, \qquad + \,
j^0 (t, (R_4 w)(t), Mw(t); Mv - Mw(t)) \ge 0
\ \ \mbox{\rm for all} \ \, v \in K.
\end{cases}\leqno{(*)}
\end{equation*}
The main feature of the problem $(*)$ is the explicit
dependence of the data $A$, $f$, $\varphi$ and $j$
on both the time parameter and the history-dependent operators.
The inequality $(*)$ without history-dependent operators is called a time-dependent variational-hemivariational inequality.
The problem $(*)$ is a history-dependent inequality
when $A$, $f$, $\varphi$ and $j$ do depend explicitly
on the operators $R_1$, $R_2$, $R_3$, and $R_4$.
On the other hand, the problem $(*)$ is called
a time-dependent variational inequality
if $j\equiv 0$ and a time-dependent
hemivariational inequality if $\varphi \equiv 0$.
It is also called quasi variational-hemivariational
inequality since the convex potential $\varphi$
depends on the solution itself.

The notion of a history-dependent operator was
introduced in~\cite{SMatei2011} and used, for instance,
in a few recent
papers~\cite{Mig2021,MOS13,MOS18,SHM2015,SMatei,SPatrulescu,SXiao2015} and the references therein.
All results available in the literature are obtained
only for some special versions of the history-dependent
variational-hemivariational inequalities. These particular
cases of $(*)$ are mainly explored in the space $C(I; K)$ of continuous functions defined on $I$ with values in $K$,
where $I$ denotes either bounded or unbounded time
interval.
For instance, when the operator $A$,
the functions $\varphi$ and $j$ are time independent,
and the operators $R_1$, $R_2$ and $R_4$ are not present,
then the problem $(*)$ has been treated in~\cite{SM2016} and~\cite[Chapter 6]{SM2}.
History-dependent variational inequalities studied
in~\cite{SMatei2011,SMatei} correspond to the inequality
$(*)$ with $j \equiv 0$ and the time independent operator $A$,
and the function $\varphi$,
and without operators $R_1$, $R_2$ and $R_4$.
The results on the inequality $(*)$
with no explicit dependence on time,
and without $j$, $R_2$ and $R_4$ can be found
in~\cite{SXiao2015}.

The inequality in problem $(*)$ has not been considered,
to the best of the authors' knowledge, in the literature
so far.
In this regard, our Theorem~\ref{theorem1q} is completely new and
has no analogues whatever.
Further, even without history-dependent operators,
we are not aware of any result for the time-parameter dependent inequality $(*)$.
The first novelty of the paper is Theorem~\ref{theorem1q}, on existence and uniqueness of solution.
The proof is based on
a recent result in~\cite{MOS30} for an elliptic problem,
and a fixed point argument in~\cite{SM2}.
The second novelty of the paper is the analysis of a new quasistatic nonsmooth frictional contact model in viscoelasticity
for which we establish existence, uniqueness and regularity
of weak solutions.
The contact boundary conditions involve a unilateral Signorini condition for the normal velocity combined with the nonmonotone normal damped response condition, the Coulomb law of dry friction,
and the normal compliance condition.
Some results on the evolutionary counterparts
of the variational-hemivariational inequality $(*)$
can be found in~\cite{hms,KULIG,MigJOTA,Mig2021,MigBai,MSZ2019,MBZ,MigZENG,SM2}.
We refer also to~\cite{HS,NP,SST,SMatei}
and the references therein for the recent results on the mathematical theory of contact mechanics and related issues.


Throughout the paper, we adopt the following notation.
Let $X$ be a Banach space with the norm
$\| \cdot \|_X$, $X^*$ denote its dual space and  $\langle\cdot,\cdot\rangle_{X^* \times X}$
be the duality brackets between $X^*$ and $X$.
Often, when no confusion arises, for simplicity, we skip the subscripts.
The symbols $``\rightarrow"$ and
$``\rightharpoonup"$ denote
the strong and the weak convergence, respectively.
A space $X$ endowed with the weak topology is denoted
by $X_w$.
Given $0 < T < \infty$ and a subset $K \subset X$,
we denote by $L^2(0, T; K)$ the set (equivalence classes) of functions in $L^2(0, T; X)$ that for almost everywhere
$t \in (0, T)$ have values in $K$.
We denote by $C(0, T; K)$ the set of continuous functions
on $[0, T]$ with values in $K$.
We write $C(0,T)$ when $K = \real$.
Let ${\mathcal L}(E, F)$ stand for
the space of linear bounded operators from
a Banach space $E$ to a Banach space $F$.
It is well known that
${\mathcal L}(E, F)$ is a Banach space endowed
with the usual norm $\| \cdot \|_{{\mathcal L}(E, F)}$.
The adjoint to $A \in {\mathcal L(E, F)}$ is denoted by
$A^* \in {\cal L}(F^*, E^*)$.

The rest of the paper is organized as follows.
In Section~\ref{s1}, we state the main result
on existence and uniqueness of solution to the
abstract quasi variational-hemivariatio\-nal inequality involving history-dependent operators.
We also provide comments and remarks on the main hypotheses.
Section~\ref{s2} is devoted to the proof of the main theorem.
Finally, in Section~\ref{Application},
we work in detail a quasistatic contact problem with history-dependent operators.

\section{Existence and uniqueness result}\label{s1}
\setcounter{equation}0

\noindent
The goal of this section is to study the existence and uniqueness of solution to the quasi
variational-hemivariatio\-nal inequality involving
history-dependent operators.

Let $E$, $X$, $Y$, and $Z$ be separable Banach
spaces, and $V$ be a separable and reflexive Banach space.
The norm in $V$ and the duality brackets
between $V^*$ and $V$ are denoted by
$\| \cdot \|$ and $\langle\cdot,\cdot\rangle$,
respectively.
\begin{Problem}\label{Problem1q}
Find $w \in L^2(0, T; V)$ such that $w(t) \in K$
for a.e. $t \in (0, T)$ and
\begin{equation*}
\begin{cases}
\displaystyle
\langle A (t, (R_1 w)(t), w(t)) - f(t, (R_2 w)(t)),
v - w(t) \rangle \\[1mm]
\qquad + \,
\varphi (t, (R_3 w)(t), w(t), v)
- \varphi(t, (R_3 w)(t), w(t), w(t)) \\[1mm]
\ \ \ \, \qquad + \,
j^0 (t, (R_4 w)(t), Mw(t); Mv - Mw(t)) \ge 0
\ \ \mbox{\rm for all} \ \, v \in K,
\ \mbox{\rm a.e.} \ t \in (0, T).
\end{cases}
\end{equation*}
\end{Problem}

We posit the following hypotheses for the data of  Problem~\ref{Problem1q}.

\smallskip
\noindent
$\underline{H(A)}:$ \quad
$\displaystyle A \colon [0, T] \times E \times V \to V^*$
is such that

\smallskip

\lista{
\item[(a)]
$A(\cdot, \cdot, v)$ is continuous for all $v \in V$,
\smallskip
\item[(b)]
$\| A(t, \lambda, v) \|_{V^*} \le a_0(t)
+ a_1 \| \lambda \|_E + a_2 \| v \|$ for all $t \in [0, T]$,
$\lambda \in E$, $v \in V$ with $a_0\in C(0, T)_+$,
 $a_1$, $a_2 \ge 0$,
\smallskip
\item[(c)]
$A(t, \lambda, \cdot)$ is demicontinuous
for all $t \in [0,T]$, $\lambda \in E$,
\smallskip
\item[(d)]
$\langle A (t, \lambda_1, v_1)- A (t, \lambda_2, v_2),
v_1 - v_2 \rangle_{V^* \times V}
\ge m_A \, \| v_1 - v_2 \|^2 - {\bar{m}}_A \,
\| \lambda_1 - \lambda_2 \|_E \| v_1 - v_2\|
$ \\[1mm]
for all $t \in [0, T]$, $\lambda_1$, $\lambda_2 \in E$,
$v_1$, $v_2 \in V$ with $m_A > 0$, ${\bar{m}}_A \ge 0$.
}

\smallskip
\noindent
$\underline{H(f)}:$ \quad
$\displaystyle f \colon [0, T] \times X \to V^*$
is such that

\smallskip

\lista{
	\item[(a)]
	$f(\cdot, \xi)$ is continuous for all $\xi \in X$,
	\smallskip
	\item[(b)]
	$\| f(t, \xi_1) - f(t, \xi_2)\|_{V^*} \le
	L_f \, \| \xi_1 - \xi_2\|_X$ for all $t \in [0, T]$,
	$\xi_1$, $\xi_2 \in X$ with $L_f > 0$.
}

\smallskip

\noindent
$\underline{H(\varphi)}:$ \quad
$\varphi \colon [0, T] \times Y \times V \times V \to \real$ is such that

\smallskip

\lista{
\item[(a)]
$\varphi(t, \eta, w, \cdot)$ is convex, lower semicontinuous  for all
$t \in [0, T]$, $\eta \in Y$, $w \in V$,
\smallskip
\item[(b)]
$
\varphi (t, \eta_1, w_1, v_2) -
\varphi (t, \eta_1, w_1, v_1) +
\varphi (t, \eta_2, w_2, v_1) -
\varphi (t, \eta_2, w_2, v_2) \\[1mm]
$
$
~\qquad \le \alpha_{\varphi} \,
\| w_1 - w_2 \| \| v_1 - v_2 \|
+ \beta_{\varphi}\| \eta_1 - \eta_2\|_Y \| v_1 - v_2\|
$ \\[1mm]
for all $t \in [0, T]$, $\eta_1$, $\eta_2 \in Y$,
$w_1$, $w_2$, $v_1$, $v_2 \in V$
with $\alpha_{\varphi}$, $\beta_{\varphi} \ge 0$.
}

\smallskip

\noindent
$\underline{H(j)}:$ \quad
$j \colon [0, T] \times Z \times X \to \real$ is such that

\smallskip

\lista{
\item[(a)]
$j(\cdot, \cdot, v)$ is continuous for all $v \in X$,  \smallskip
\item[(b)]
$j(t, \zeta, \cdot)$ is locally Lipschitz
for all $t \in [0, T]$, $\zeta \in Z$,
\smallskip
\item[(c)]
$\| \partial j(t, \zeta, v) \|_{X^*} \le
c_{0j}(t) + c_{1j} \| \zeta \|_Z + c_{2j} \| v \|_X$
\smallskip
for all $t \in [0, T]$, $\zeta \in Z$, $v \in X$
with
$c_{0j} \in C(0, T)_+$,  $c_{1j}$, $c_{2j} \ge 0$, \smallskip
\item[(d)]
$j^0(t, \zeta_1, v_1; v_2 - v_1)
+ j^0(t, \zeta_2, v_2; v_1 - v_2)
\le m_j \, \| v_1 - v_2 \|^2_X + m_1 \| \zeta_1 - \zeta_2 \|_Z
\| v_1 - v_2 \|_X
$ \\[1mm]
for all $t \in [0, T]$, $\zeta_1$, $\zeta_2 \in Z$,
$v_1$, $v_2 \in X$ with $m_j$, $m_1 \ge 0$.
}

\smallskip

\noindent
$\underline{H(K)}:$ \quad
$K$ is a nonempty, closed and convex subset of $V$.

\smallskip

\noindent
$\underline{H(M)}:$ \quad
$M \colon V \to X$
is a linear, bounded and compact operator.

\smallskip

\noindent
$\underline{H(R)}:$ \quad
$R_1 \colon L^2(0, T; V) \to L^2(0, T; E)$,
$R_2 \colon L^2(0, T; V) \to L^2(0, T; X)$,

\qquad \ \
$R_3 \colon L^2(0, T; V) \to L^2(0, T; Y)$
and
$R_4 \colon L^2(0, T; V) \to L^2(0, T; Z)$
are such that


\lista{
	\item[(a)]
	$\displaystyle
	\| (R_1 v_1)(t) - (R_1 v_2)(t) \|_{E} \le
	c_{R_1} \int_0^t \| v_1(s) - v_2(s) \| \, ds$ \\ [1mm]
	for all $v_1$, $v_2 \in L^2(0, T; V)$,
	a.e. $t\in (0, T)$ with $c_{R_1} > 0$, \smallskip
	\item[(b)]
	$\displaystyle
	\| ({R_2} v_1)(t) - ({R_2} v_2)(t) \|_X \le
	c_{R_2} \int_0^t \| v_1(s) - v_2(s) \| \, ds$ \\ [1mm]
	for all $v_1$, $v_2 \in L^2(0, T; V)$,
	a.e. $t\in (0, T)$ with $c_{R_2} > 0$, \smallskip
	\item[(c)]
	$\displaystyle
	\| ({R_3} v_1)(t) - ({R_3} v_2)(t) \|_Y \le
	c_{R_3} \int_0^t \| v_1(s) - v_2(s) \| \, ds$
	\\ [1mm]
	for all $v_1$, $v_2 \in L^2(0, T; V)$,
	a.e.\ $t\in (0, T)$ with $c_{R_3} > 0$,
	\item[(d)]
	$\displaystyle
	\| ({R_4} v_1)(t) - ({R_4} v_2)(t) \|_Z \le
	c_{R_4} \int_0^t \| v_1(s) - v_2(s) \| \, ds$
	\\ [1mm]
	for all $v_1$, $v_2 \in L^2(0, T; V)$,
	a.e.\ $t\in (0, T)$ with $c_{R_4} > 0$.
}

\smallskip

\noindent
$\underline{(H_0)}:$ \quad
$m_{j} \| M \|^2 + \alpha_{\varphi}< m_A$.

\medskip

Before we state the main result of the paper, we
discuss the hypotheses.
As concerns hypothesis $H(A)$(c),
an operator $A \colon V \to V^*$
is said to be demicontinuous if it is strongly-weakly continuous.
In hypothesis $H(j)$(c), (d), and in what follows,
the symbols $\partial j$ and $j^0$ stand for
the generalized subgradient and the generalized directional derivative of a function $j$ with respect to the last variable.
We recall,
see~\cite{cla,DMP1,MOSBook},
that given a locally Lipschitz function
$\psi \colon X \to \real$ on a Banach space $X$,
the generalized subgradient of $\psi$ at $x$,
is given by
\begin{equation*}
\partial \psi (x) = \{\, x^* \in X^* \mid
{\langle x^*, v \rangle} \le \psi^{0}(x; v)
\ \ \mbox{for all} \ \ v \in X \, \},
\end{equation*}
where the generalized directional derivative of $\psi$ at
$x \in X$ in the direction $v \in X$
is defined by
\begin{equation*}
\psi^{0}(x; v) = \limsup_{y \to x, \ \lambda \downarrow 0}
\frac{\psi(y + \lambda v) - \psi(y)}{\lambda}.
\end{equation*}
Further, we use the notation
$\| {\mathscr S} \|_X
= \sup \{ \, \| u \|_X \mid u \in {\mathscr S} \, \}$
for any set ${\mathscr S} \subset X$,

The operators that satisfy inequalities in $H(R)$
are called history-dependent operators since their current value for a given function at the time instant $t$
depends on the values of the function at the moments
from the time interval $[0,t]$. Such operators include
Volterra-type operators and other integral-type operators,
see~\cite{SM2} and the literature therein.
There are four history-dependent operators in the inequality
which are presented in the operator $A$, the function $f$,
and two potentials $\varphi$ and $j$, respectively.

The function $\varphi$, in Problem~\ref{Problem1q}, is convex in the last variable and it provides a variational term, while the function $j$ is locally Lipschitz in the last argument and
it introduces a hemivariational part. Since the convex potential
$\varphi$ depends on the solution in the third argument,
the inequality is called quasi variational-hemivariational
inequality.
We note that the condition $w(t) \in K$ for a.e. $t \in (0, T)$
imposes an additional constraint in the problem.

\begin{Remark}\label{REM2}
The following conditions can be used to check
the hypothesis $H(A)${\rm (d)}.
If $A(t, \lambda, \cdot)$
is strongly monotone with $m_A>0$
and $A(t,\cdot, v)$ is
Lipschitz with $L_A>0$, i.e.,
\begin{eqnarray*}
&&
\langle A (t, \lambda, v_1)
- A(t, \lambda, v_2), v_1 - v_2 \rangle
\ge m_A \| v_1 - v_2 \|^2
\ \ \mbox{for all} \ \  \lambda \in E,
\ v_1, v_2 \in V, \\[2mm]
&&
\| A (t, \lambda_1, v)
- A(t, \lambda_2, v) \|_{V^*} \le L_A\,
\| \lambda_1 - \lambda_2\|_E
\ \ \mbox{for all} \ \
\lambda_1, \lambda_2 \in E,
\ v \in V,
\end{eqnarray*}
for a.e. $t \in (0,T)$,
then
\begin{equation*}
\langle A (t, \lambda_1, v_1)
- A(t, \lambda_2, v_2), v_1 - v_2 \rangle
\ge m_A \, \| v_1 - v_2 \|^2
- L_A \, \|\lambda_1 - \lambda_2 \|_E
\| v_1 - v_2\|
\end{equation*}
for all $\lambda_1$, $\lambda_2 \in E$,
$v_1$, $v_2 \in V$, a.e. $t \in (0, T)$.
\end{Remark}

We also need stronger assumptions.

\smallskip

\smallskip

\noindent
$\underline{H(\varphi)_1}:$ \quad
$\varphi \colon [0, T] \times Y \times V \times V \to \real$ satisfies $H(\varphi)$(a),(b) and

\smallskip

\lista{
\item[(c)]
$\varphi (t, \eta, w, v_1) - \varphi(t, \eta, w, v_2)
\le
\left(
c_{\varphi_1}(t) + c_{\varphi_2}(\| w \|) + c_3 \| \eta \|_Y
\right) \| v_1 - v_2\|$ \\[2mm]
for all $(t, \eta, w) \in [0,T] \times Y \times V$,
$v_1$, $v_2 \in V$, where
$c_{\varphi_1} \colon [0, T] \to [0,\infty)$ and
$c_{\varphi_2} \colon [0, \infty) \to [0, \infty)$
are continuous functions, and $c_3 > 0$.
\smallskip
\item[(d)]
$
\limsup \left(
\varphi (t_n, \eta_n, w_n, v) -
\varphi (t_n, \eta_n, w_n, w_n) \right)
\le
\varphi (t, \eta, w, v) - \varphi (t, \eta, w, w)
$ \\[2mm]
for any $v \in V$,
$t_n \to t$ in $[0, T]$,
$\eta_n \to \eta$ in $Y$ and
$w_n \rightharpoonup w$ in $V$.
}

\smallskip

\noindent
$\underline{H(j)_1}:$ \quad
$j \colon [0, T] \times Z \times X \to \real$
satisfies $H(j)$ (a)--(d) and
\begin{equation}\label{EE1}
\limsup j^0(t_n, \zeta_n, Mv; Mv-Mw_n)
\le j^0(t, \zeta, Mv; Mv-Mw)
\end{equation}
for any $v \in V$,
$t_n \to t$ in $[0, T]$,
$\zeta_n \to \zeta$ in $Z$ and
$w_n \rightharpoonup w$ in $V$.

\smallskip

{\rm
We comment on hypotheses $H(\varphi)_1$.
Consider a function $\varphi$ which is independent
of the first three variables, that is,
$\varphi(t, \eta, w, v) = \varphi (v)$.
Then, under $H(\varphi)$(a),
hypotheses $H(\varphi)$(b) and (d)
can be omitted.
It is clear that in this case, condition (b)
is trivially satisfied with
$\alpha_\varphi = \beta_{\varphi} = 0$.
The condition $H(\varphi)$(d) is a consequence of $H(\varphi)$(a).
Indeed,
for any $v \in V$ and
$w_n \rightharpoonup w$ in $V$,
by the weak lower semicontinuity of $\varphi$, we have
$\varphi(w) \le \liminf \varphi(w_n)$.
Hence
$$
\limsup \big(
\varphi (v) - \varphi(w_n) \big)
= \varphi(v) + \limsup \, (-\varphi(w_n))
\le \varphi(v) - \varphi(w).
$$
Remark also that if
$\varphi(t, \eta, w, v) = \varphi(w, v)$,
then condition $H(\varphi)$(b)
was already used in~\cite{MOS30,SM2}
and the references therein,
conditons $H(\varphi)$(b) and (c) together were used
in~\cite[Theorem~10]{Xiao2019},
while a version of $H(\varphi)$(d)
was employed in~\cite{MOTSOF,Sofonea}, respectively.
}

Next, we provide a prototype of a function $j$
that satisfies $H(j)_1$.
\begin{Remark}\label{RR1}
Let $j \colon [0, T] \times Z \times X \to \real$
be given by	
$$
j(t, \zeta, v) = \alpha(t, \zeta) \, g(v)
\ \ \mbox{for} \ \
(t, \zeta, v) \in [0, T] \times Z \times X,
$$
where

\smallskip

\noindent
$\underline{H(\alpha)}:$ \quad
$\alpha \colon [0, T] \times Z \to [0, \infty)$
is a continuous function such that

\smallskip

\qquad \quad
$\alpha(t, \zeta) \le \alpha_0$ for
$(t, \zeta) \in [0, T] \times Z$ with $\alpha_0 > 0$.

\smallskip

\noindent
$\underline{H(g)}:$ \quad
$g \colon X \to \real$ is a locally Lipschitz
function such that

\smallskip

\lista{
\item[\rm (a)]
$\| \partial g (v) \|_X \le c_{0g} + c_{1g} \| v \|_X$
for
$v \in X \ \mbox{with} \ c_{0g}, c_{1g} \ge 0$,
\smallskip
\item[\rm (b)]
$g^0(v_1; v_2-v_1) + g^0(v_2; v_1-v_2) \le
m_g \|v_1 - v_2\|^2_X$ for $v_1$, $v_2 \in X$ with
$m_g \ge 0$.
}

\smallskip

\noindent
Under $H(M)$ the function $j$ satisfies $H(j)_1$.
The hypotheses $H(j)$ are obviously satisfied.
We show $(\ref{EE1})$. Indeed, let $v \in V$,
$t_n \to t$ in $\real$,
$\zeta_n \to \zeta$ in $Z$ and
$w_n \rightharpoonup w$ in $V$.
By~\cite[Proposition~3.23(ii)]{MOSBook},
we have
\begin{equation}\label{EE2}
\limsup g^0(Mv, Mv-Mw_n) \le g^0(Mv; Mv-Mw).
\end{equation}
Then
\begin{eqnarray*}
&&\hspace{-0.5cm}
j^0(t_n, \zeta_n, Mv; Mv-Mw_n)
= \alpha(t_n, \zeta_n) \, g^0(Mv; Mv-Mw_n) \\[2mm]
&&
\le |\alpha(t_n, \zeta_n) - \alpha(t, \zeta)|
(c_{0g} + c_{1g} \| v \|_X) \| M \| \| v - w_n \|
+ \alpha (t, \zeta) \, g^0(Mv; Mv-Mw_n)
\end{eqnarray*}
which combined with $(\ref{EE2})$ implies $(\ref{EE1})$.	
\end{Remark}

Further, we note that condition $H(g)$(b) is equivalent
to the condition
\begin{equation}\label{RELAXED}
\langle \partial g(v_1) - \partial g(v_2), v_1-v_2 \rangle
\ge - m_g \| v_1-v_2\|^2_X
\ \ \mbox{for all} \ \ v_1, v_2 \in X.
\end{equation}
The latter is known as the relaxed monotonicity condition,
see~\cite{MOSBook} and the references therein. Examples
of nonconvex functions which satisfy $H(g)$ can be found
in~\cite{MOS18}. If $g \colon X \to \real$ is a convex
function, then $H(g)$(b) and (\ref{RELAXED}) hold
with $m_g = 0$ thanks to the monotonocity of the convex subdifferential.

\medskip

In the proof of the main result we need the
following observation.
\begin{Lemma}\label{LL1}
Under hypotheses $H(A)$, $H(j)$, $H(M)$ and $(H_0)$,
the multivalued map
$$
A(t, \lambda, \cdot) + M^* \partial j (t, \zeta, M \cdot)
\colon V \to 2^{V^*} \setminus \{ \emptyset \}
$$
is monotone for all $(t, \lambda, \zeta) \in [0, T] \times E \times Z$.
\end{Lemma}

\noindent
{\bf Proof.} \
Let $(t, \lambda, \zeta) \in [0, T] \times E \times Z$
be fixed,
$v_1$, $v_2 \in V$, and
$z_1 \in \partial j(t, \zeta, M v_1)$,
$z_2 \in \partial j(t, \zeta, M v_2)$.
By the definition of the generalized subgradient
and $H(j)$, we have
\begin{eqnarray*}
&&\hspace{-0.4cm}
\langle M^* z_1 - M^* z_2, v_2 - v_1 \rangle
=
\langle z_1, M (v_2-v_1) \rangle_{X} +
\langle z_2, M (v_1-v_2) \rangle_{X} \\[2mm]
&&\quad
\le j^0(t, \zeta, Mv_1; Mv_2 - Mv_1)
+ j^0(t, \zeta, Mv_2; Mv_1 - Mv_2)
\le m_j \| M \|^2 \| v_1 - v_2 \|^2.
\end{eqnarray*}
Hence, by $(H_0)$, we get
\begin{eqnarray*}
&&\hspace{-0.3cm}
\langle A(t, \lambda, v_1) + M^* z_1
- A(t, \lambda, v_2) - M^* z_2, v_1 - v_2 \rangle
= \langle A(t, \lambda, v_1) - A(t, \lambda, v_2),
v_1-v_2 \rangle \\[2mm]
&& \quad
+
\langle M^* z_1 - M^* z_2, v_1 - v_2 \rangle
\ge (m_A - m_j \| M \|^2) \| v_1 - v_2\|^2 \ge 0,
\end{eqnarray*}
which completes the proof.
\hfill$\Box$

\medskip

We conclude this section with the main
existence and uniqueness result of the paper.
Its proof will be given in the next section.
\begin{Theorem}\label{theorem1q}
	Under hypotheses $H(A)$, $H(f)$, $H(\varphi)_1$,
	$H(j)_1$, $H(K)$, $H(M)$, $H(R)$, and $(H_0)$,
	Problem~{\rm \ref{Problem1q}} has a unique solution
	$w \in L^2(0, T; V)$ such that $w(t) \in K$
	for a.e. $t \in (0, T)$.
\end{Theorem}

\section{Proof of the main result}\label{s2}

The proof of Theorem~\ref{theorem1q}
will be performed in six steps.

\medskip

\noindent
{\bf Step 1}.
(Unique solvability of an auxiliary problem).

We shall prove that under hypotheses
$H(A)$, $H(f)$, $H(\varphi)$, $H(j)$, $H(K)$, $H(M)$, and $(H_0)$,
for any fixed
$(t, \lambda, \xi, \eta, \zeta) \in
[0, T] \times E \times X \times Y \times Z$,
the following auxiliary problem:
find $w \in K$ such that
\begin{equation}\label{PP4}
\begin{cases}
\displaystyle
\langle A (t, \lambda, w) - f(t,\xi), v - w \rangle
+ \, \varphi (t, \eta, w, v)
- \varphi(t, \eta, w, w) \\[1mm]
\ \ \ \qquad + \,
j^0 (t, \zeta, Mw; Mv - w) \ge 0
\ \ \mbox{\rm for all} \ \, v \in K
\end{cases}
\end{equation}
has a unique solution $w \in K$.

We will apply~\cite[Theorem~18]{MOS30} established
for an elliptic variational-hemivariational inequality.
We define
${\widetilde A} \colon V \to V^*$,
${\widetilde \varphi} \colon V \times V \to \real$,
${\widetilde j} \colon V \to \real$, and
${\widetilde f} \in L^2(0, T; V^*)$ by
\begin{equation*}
{\widetilde A} v = A(t, \lambda, v), \ \
{\widetilde \varphi}(z, v) = \varphi (t, \eta, z, v),
\ \
{\widetilde j}(v) = j(t, \zeta, Mv),
\ \
{\widetilde f}(t) = f(t, \xi)
\end{equation*}
for $t \in [0, T]$, $z$, $v \in V$.
We will verify that the above data satisfy conditions (22)--(26), (29) and (30) of~\cite[Theorem~18]{MOS30}.
It is easy to check, by $H(A)$, that
$$
\langle {\widetilde A} v_1 - {\widetilde A} v_2, v_1 - v_2
\rangle =
\langle A(t, \lambda, v_1) - A(t, \lambda, v_2), v_1 - v_2
\rangle \ge m_A \, \| v_1 - v_2 \|^2
$$
which means that ${\widetilde A}$ is strongly monotone with
constant $m_A > 0$.
Using~\cite[Remark~13, p.147]{SM2}, it follows that
$$
\langle {\widetilde A} v, v - v_0\rangle \ge
m_A \| v \|^2 - \beta_1 \| v\| - \beta_2
\ \ \mbox{for all} \ \ v \in V, \ \mbox{any} \
v_0 \in K \ \mbox{with} \ \beta_1, \beta_2 \in \real.
$$
From $H(A)$(b) and (c), it is clear that ${\widetilde A}$
is bounded and demicontinuous, so also pseudomonotone,
see, e.g.,~\cite[Theorem~3.69(i)]{MOSBook}.
Therefore, condition~(22) in~\cite{MOS30} holds.
By $H(\varphi)$(a) and (b) we easily
deduce condition~(23) in~\cite{MOS30} with
$\alpha_{\varphi} > 0$.

Hypothesis $H(j)$ combined with the chain rule for the generalized directional derivative, see~\cite[Proposition~3.37(i)]{MOSBook},
implies that condition~(24) in~\cite{MOS30}
is satisfied with $\alpha_j = m_j \| M \|^2$.
The conditions (25) and (26) in~\cite{MOS30}
are consequences of $H(K)$ and $H(f)$, respectively.
Finally, the smallness conditions (29) and (30)
in~\cite{MOS30} hold due to $(H_0)$.
Having verified the above hypotheses,
we are now in a position to apply~\cite[Theorem~18]{MOS30}
to deduce that problem $(\ref{PP4})$
has a unique solution $w \in K$.

\medskip

\noindent
{\bf Step 2}.
(A priori estimate for the auxiliary problem).

Let $(t, \lambda, \xi, \eta, \zeta) \in
[0, T] \times E \times X \times Y \times Z$
be fixed.
In addition to hypotheses in Step~1, we suppose that $H(\varphi)_1$$(c)$ holds.
We shall demonstrate that if $w \in K$ solves the auxiliary
problem $(\ref{PP4})$, then
\begin{equation}\label{EST33}
\| w \| \le C \left(1+
\| \lambda \|_E + \| \xi \|_X + \|\eta\|_Y +
\| \zeta \|_Z \right)
\ \ \mbox{with a constant} \ \ C > 0.
\end{equation}

For the proof of (\ref{EST33}), let $v_0$ be any element in $K$.
We choose $v = v_0$ in the inequality (\ref{PP4}) to get
\begin{equation*}
\langle A (t, \lambda, w) - f(t,\xi), v_0 - w \rangle
+ \, \varphi (t, \eta, w, v_0)
- \varphi(t, \eta, w, w) +
j^0 (t, \zeta, Mw; Mv_0 - Mw) \ge 0.
\end{equation*}
We shall  estimate separately each term in this inequality.
From $H(A)$, we obtain
\begin{eqnarray}
&&\langle A(t, \lambda, w), w - v_0 \rangle
=
\langle A(t, \lambda, w) - A(t, \lambda, v_0), w - v_0
\rangle +
\langle A(t, \lambda, v_0), w - v_0 \rangle
\nonumber \\[2mm]
&&\quad
\ge m_A \| w - v_0 \|^2 +
\langle A(t, \lambda, v_0), w - v_0 \rangle. \label{AA1}
\end{eqnarray}
Next, we conclude from $H(j)$(c) and (d) that
\begin{eqnarray}
&&\hspace{-0.6cm}
j^0(t, \zeta, Mw; Mv_0 - Mw) =
j^0(t, \zeta, Mw; Mv_0 - Mw) +
j^0(t, \zeta, Mv_0; Mw - Mv_0) \nonumber
\\[2mm]
&& \hspace{-0.4cm}
- \,
j^0(t, \zeta, Mv_0; Mw - Mv_0) \le m_j \| M w - Mv_0 \|^2_X
\nonumber \\[2mm]
&&\hspace{-0.4cm} \quad
+ \left(
c_{0j}(t) + c_{1j} \| \zeta \|_Z + c_{2j} \| Mv_0 \|_X
\right) \, \| Mw - Mv_0\|_X \nonumber \\[2mm]
&& \hspace{-0.4cm}
\qquad
\le m_j \| M \|^2 \| w - v_0 \|^2
+ \| M \| \left(
c_{0j}(t) + c_{1j} \| \zeta \|_Z + c_{2j} \| M \| \| v_0 \|
\right) \, \| w-v_0\|. \label{jj1}
\end{eqnarray}
We use hypotheses $H(\varphi)_1$(b) and (c)
and the triangle
inequality to get
\begin{eqnarray}
&&\hspace{-0.8cm}
\varphi (t, \eta, w, v_0) -
\varphi (t, \eta, w, w) \le
\varphi (t, \eta, z_0, v_0) -
\varphi (t, \eta, z_0, w)
+ \alpha_{\varphi} \| w - z_0\| \| w-v_0\|
\nonumber
\\[2mm]
&& \hspace{-0.8cm}
\quad
\le
(c_{\varphi_1}(t) + c_{\varphi_2}(\| z_0 \|)
+ c_3 \| \eta \|_Y) \| w - v_0\| +
\alpha_{\varphi} \left( \| w-v_0\| + \| v_0 - z_0\| \right) \,
\| w - v_0\|
\nonumber \\[2mm]
&&\hspace{-0.8cm}\qquad
\le \alpha_{\varphi} \| w - v_0 \|^2
+
(c_{\varphi_1}(t) + c_{\varphi_2}(\| z_0 \|)
+ c_3 \| \eta \|_Y + \alpha_{\varphi} \| v_0 - z_0\|)
\| w-v_0 \| \label{fifi1}
\end{eqnarray}
for any $z_0 \in V$.
The hypotheses $H(A)$ and $H(f)$ imply that
\begin{eqnarray}
&&
\| A(t, \lambda, v_0) \|_{V^*} \le
a_0(t) + a_1 \| \lambda\|_E + a_2 \| v_0\|,
\label{AAA}
\\[2mm]
&&
\| f(t, \xi) \|_{V^*} \le
L_f \, \| \xi \|_X + L_f \, \|\xi_0 \|_X
+ \|f(t, \xi_0)\|_{V^*}\label{ff1}
\end{eqnarray}
for any element $\xi_0 \in X$.
We combine the estimates (\ref{AA1})--(\ref{ff1})
to obtain
\begin{eqnarray*}
&&
(m_A - m_j \| M \|^2 - \alpha_{\varphi}) \| w - v_0 \|
\le
a_0(t) + \|f(t, \xi_0)\|_{V^*} + \| M \| c_{0j}(t)
\\[2mm]
&&\quad
+ \,
c_{\varphi_1}(t) + c_{\varphi_2}(\| z_0 \|)
+ c_3 \| \eta \|_Y + a_2 \| v_0\| + L_f \|\xi_0\|_X
+ c_{2j} \| M \|^2 \| v_0\| \\[2mm]
&&\qquad
+ \, \alpha_{\varphi} \| v_0 - z_0\|
+ a_1 \| \lambda\|_E + L_f \|\xi \|_X +
c_{1j} \| M \| \|\zeta\|_Z
\label{EST2}
\end{eqnarray*}
for any $v_0 \in K$, $z_0 \in V$ and $\xi_0 \in X$.
Hence, there is a constant $C > 0$ such that
$$
(m_A - m_j \| M \|^2 - \alpha_{\varphi}) \| w - v_0 \|
\le
C\left(
1 + \gamma(t) + \|\lambda\|_E + \|\xi \|_X + \| \eta\|_Y + \| \zeta\|_Z \right),
$$
where
$$
\gamma(t) = a_0(t) + \| f(t, \xi_0) \|_X +
\| M \| c_{0j}(t) + c_{\varphi_1}(t).
$$
We have
$\gamma \in C(0, T)$ and
$\gamma(t) \le \gamma_0$ for all $t \in [0, T]$
with $\gamma_0 > 0$.
Finally, by the smallness condition $(H_0)$,
we deduce the bound (\ref{EST33}).

\medskip

\noindent
{\bf Step 3}.
(The Minty formulation of the auxiliary problem).


We shall prove that under hypotheses $H(A)$, $H(f)$ $H(\varphi)$, $H(j)$, $H(K)$, $H(M)$ and $(H_0)$,
for any fixed
$(t, \lambda, \xi, \eta, \zeta) \in [0, T] \times E \times Z$,
the problem $(\ref{PP4})$ is equivalent
to the following Minty inequality:
%
%
find $w \in K$ such that
\begin{equation}\label{PP2}
\begin{cases}
\displaystyle
\langle A (t, \lambda, v) - f(t,\xi), v - w \rangle
+ \, \varphi (t, \eta, w, v)
- \varphi(t, \eta, w, w) \\[1mm]
\ \ \ \qquad + \,
j^0 (t, \zeta, Mv; Mv - Mw) \ge 0
\ \ \mbox{\rm for all} \ \, v \in K.
\end{cases}
\end{equation}

Indeed, let $w \in K$ be a solution to (\ref{PP4})
and $v \in K$.
We apply Lemma~\ref{LL1} and
for all
$\eta_v \in \partial j(t, \zeta, Mv)$ and
$\eta_w \in \partial j(t, \zeta, Mw)$, we compute
\begin{eqnarray*}
	&&\hspace{-0.7cm}
	\langle A (t, \lambda, v) - f(t,\xi), v - w \rangle
	+ \, \varphi (t, \eta, w, v) - \varphi(t, \eta, w, w)
	+ \, j^0 (t, \zeta, Mv; Mv - Mw) \\[1mm]
	&&
	\ge  \langle A (t, \lambda, v) - A (t, \lambda, w), v - w \rangle
	+ \langle A (t, \lambda, w) - f(t,\xi), v - w \rangle
	\\[1mm]
	&&
	+ \, \varphi (t, \eta, w, v) - \varphi(t, \eta, w, w)
	+ \langle \eta_v-\eta_w, Mv-Mw \rangle
	+ \, j^0 (t, \zeta, Mw; Mv - Mw) \\[1mm]
	&&
	= \langle A (t, \lambda, v) + M^* \eta_v
	- A (t, \lambda, w) - M^* \eta_w, v - w \rangle
	+ \langle A (t, \lambda, w) - f(t,\xi), v - w \rangle
	\\[1mm]
	&&
	+ \, \varphi (t, \eta, w, v) - \varphi(t, \eta, w, w)
	+ \, j^0 (t, \zeta, Mw; Mv - Mw) \ge 0.
\end{eqnarray*}
Thus $w \in K$ solves (\ref{PP2}).
Conversely, let $w \in K$ solve (\ref{PP2}).
Let $z \in K$, $\theta \in (0, 1)$ and
$v_\theta = \theta z + (1-\theta)w$. Then
$v_\theta \in K$ and $v_\theta = w +\theta (z-w)$.
Next, we take $v = v_\theta$ in (\ref{PP2}) to obtain
\begin{equation}\label{PP3}
\begin{cases}
\displaystyle
\langle A (t, \lambda, v_\theta) - f(t,\xi),
\theta (z-w) \rangle
+ \, \varphi (t, \eta, w, v_\theta)
- \varphi(t, \eta, w, w) \\[1mm]
\ \ \ \qquad + \,
j^0 (t, \zeta, Mv_\theta; \theta M (z - w) \ge 0.
\end{cases}
\end{equation}
From the convexity of $\varphi$ in the last variable,
the following inequality can be drawn
$$
\varphi (t, \eta, w, w+ \theta(z-w))
- \varphi(t, \eta, w, w) \le
\theta \left(
\varphi (t, \eta, w, z) - \varphi(t, \eta, w, w)
\right) .
$$
We use also the positive homogeneity of $j^0$ in its last variable, see~\cite[Proposition~3.23(i)]{MOSBook} to get
\begin{equation}\label{PP4a}
\begin{cases}
\displaystyle
\langle A (t, \lambda, v_\theta) - f(t,\xi),
z-w \rangle
+ \, \varphi (t, \eta, w, z)
- \varphi(t, \eta, w, w) \\[1mm]
\ \ \ \qquad + \,
j^0 (t, \zeta, Mv_\theta; M z - Mw) \ge 0.
\end{cases}
\end{equation}
Exploiting the compactness of $M$,
the demicontinuity of $A(t, \lambda, \cdot)$,
the upper semicontinuity of
$j^0(t, \zeta, \cdot; Mz-Mw)$,  see~\cite[Proposition~3.23(ii)]{MOSBook},
%
and passing to the limit as $\theta \to 0^+$,
we deduce
\begin{eqnarray*}
	&&
	0 \le \lim\,
	\langle A (t, \lambda, v_\theta) - f(t,\xi),
	z-w \rangle
	+ \, \varphi (t, \eta, w, z) - \varphi(t, \eta, w, w)
	\\[2mm]
	&& \quad
	+ \limsup j^0 (t, \zeta, M v_\theta; M z - Mw)
	\le
	\langle A (t, \lambda, w) - f(t,\xi),
	z-w \rangle \\[2mm]
	&&\qquad
	+ \, \varphi (t, \eta, w, z) - \varphi(t, \eta, w, w)
	+ j^0 (t, \zeta, M w; M z - Mw).
\end{eqnarray*}
Since $z \in K$ is arbitrary, we get that
$w \in K$ is a solution to problem (\ref{PP4}).
This completes the proof of Step~3.

\medskip

\noindent
{\bf Step 4}.
(Continuity of the solution map of problem (\ref{PP4})).

We shall prove that under hypotheses
$H(A)$, $H(f)$, $H(\varphi)_1$, $H(j)_1$, $H(K)$,
$H(M)$ and $(H_0)$, the solution map of (\ref{PP4})
defined by
$$
\begin{cases}
p \colon [0, T] \times E \times X \times Y \times Z
\to K, \\
p(t, \lambda, \xi, \eta, \zeta) = w
\ \ \mbox{for} \ \
(t, \lambda, \xi, \eta, \zeta)
\in [0, T] \times E \times X \times Y \times Z
\end{cases}
$$
is continuous from
$[0, T] \times E \times X \times Y \times Z$ to $V_w$.

By Step~1, the solution map for (\ref{PP4}) is well-defined and single valued.
Consider $(t_n, \lambda_n, \xi_n, \eta_n, \zeta_n) \in
[0, T] \times E \times X \times Y \times Z$,
$(t_n, \lambda_n, \xi_n, \eta_n, \zeta_n) \to
(t, \lambda, \xi, \eta, \zeta)$ in
$[0, T] \times E \times X \times Y \times Z$.
Let $w_n = p(t_n, \lambda_n, \xi_n, \eta_n, \zeta_n) \in K$
be the unique solution to problem (\ref{PP4}).
By the estimate (\ref{EST33}) of Step~2, the sequence
$\{ w_n \}$ remains in a bounded subset of $V$.
So, by the reflexivity of the latter, we may suppose
that $w_n \rightharpoonup w$ in $V$ at least for a subsequence.
Since $w_n \in K$ and, by $H(K)$, $K$ is weakly closed,
we have $w \in K$.
The Minty formulation is now useful to pass to the limit.
By Step~3, we know that
$w_n \in K$ satisfies
\begin{eqnarray}\label{Limit1}
&&
\langle A (t_n, \lambda_n, v) - f(t_n,\xi_n), v - w_n \rangle
+ \, \varphi (t_n, \eta_n, w_n, v)
- \varphi(t_n, \eta_n, w_n, w_n) \nonumber \\[1mm]
&&\ \ \ \qquad + \,
j^0 (t_n, \zeta_n, Mv; Mv - Mw_n) \ge 0
\ \ \mbox{\rm for all} \ \, v \in K.
\label{Limit1}
\end{eqnarray}
Let $v \in K$. From $H(A)$(a) and $H(f)$, we see that
\begin{eqnarray*}
&&
A(t_n, \lambda_n, v) \to A(t, \lambda, v)
\ \ \mbox{in} \ \ V^*,  \\[2mm]
&&
\| f(t_n, \xi_n)  - f(t, \xi)\|_{V^*}
\le L_f \, \| \xi_n - \xi \|_X +
\| f(t_n, \xi) - f(t, \xi) \|_{V^*} \ \to \ 0,
\end{eqnarray*}
which implies that
\begin{equation}\label{Limit2}
\langle A (t_n, \lambda_n, v) - f(t_n,\xi_n),
v - w_n \rangle \to \langle A (t, \lambda, v) - f(t,\xi),
v - w \rangle.
\end{equation}
According to $H(\varphi)_1$(c) and $H(j)_1$, we get
\begin{eqnarray}
&&\hspace{-1.0cm}
\limsup \left(
\varphi (t_n, \eta_n, w_n, v) -
\varphi (t_n, \eta_n, w_n, w_n) \right)
\le \varphi (t, \eta, w, v) - \varphi (t, \eta, w, w),
\label{Limit4} \\[2mm]
&&\hspace{-1.0cm}
\limsup j^0(t_n, \zeta_n, Mv; Mv-Mw_n)
\le j^0(t, \zeta, Mv; Mv-Mw). \label{Limit3}
\end{eqnarray}
Combining (\ref{Limit2})--(\ref{Limit4}) with inequality
(\ref{Limit1}), we deduce that
\begin{eqnarray}\label{Limit1bb}
&&
\langle A (t, \lambda, v) - f(t,\xi), v - w \rangle
+ \, \varphi (t, \eta, w, v)
- \varphi(t, \eta, w, w) \nonumber \\[1mm]
&&\ \ \ \qquad + \,
j^0 (t, \zeta, Mv; Mv - Mw) \ge 0
\ \ \mbox{\rm for all} \ \, v \in K.
\label{Limit5}
\end{eqnarray}
We use again the Minty formulation and obtain
that $w \in K$ is a solution to problem (\ref{PP4}).
From the uniqueness of solution to problem (\ref{PP4}),
we know that the whole sequence $\{ w_n \}$ converges
weakly in $V$. Thus the solution map $p$ is continuous in the aforementioned topologies.

\medskip

\noindent
{\bf Step 5}.
(The measurability of the solution map
of problem (\ref{PP4})).

We shall prove that the solution map $p$ for problem
(\ref{PP4}) defined in Step~4 is
$\Sigma \otimes {\cal B}({\mathbb X})$ measurable, where
${\mathbb X} = E \times X \times Y \times Z$,
$\Sigma$ denotes the $\sigma$-algebra of all Lebesgue measurable subsets of $[0, T]$ and
${\cal B}({\mathbb X})$ is the Borel $\sigma$-filed of ${\mathbb X}$.

In what follows, for simplicity of notation,
we will write
${\mathtt x} = (\lambda, \xi, \eta, \zeta)$ with
${\mathtt x} \in {\mathbb X}$.
Since $[0, T] \times {\mathbb X}$ is a measurbale space and
$V$ is a separable metric space,
the measurability of $p$ can be proved by an equivalent
definition of measurability by means of the distance function,
see~\cite[Proposition~4.2.4]{DMP1}.
It is enough to show that the real-valued function
\begin{equation}\label{MEAS1}
[0, T] \times {\mathbb X} \ni (t, {\mathtt x})
\ \mapsto \ \| v - p(t, {\mathtt x}) \| \in [0, \infty)
\end{equation}
is $\Sigma \otimes {\cal B}({\mathbb X})$-measurable for all $v \in V$.
The measurability of (\ref{MEAS1}) can be ensured when we prove that the set $M_r (v)$ is measurable
for all $r > 0$ and $v \in V$, where
$$
M_r(v) = \{(t, {\mathtt x}) \in [0, T] \times
{\mathbb X} \mid \| v - p(t, {\mathtt x}) \| \le r \}.
$$
Since every closed subset of
$[0, T] \times {\mathbb X}$ is measurable, it is enough to
show that $M_r(v)$ is closed for all $r > 0$ and $v \in V$.
To this end,
let $r > 0$ and $v \in V$ be fixed.
Let $(t_n, {\mathtt x}_n) \in M_r(v)$ and
$(t_n, {\mathtt x}_n) \to (t_0, {\mathtt x}_0) $
in $[0, T] \times {\mathbb X}$.
We prove that $(t_0, {\mathtt x}_0) \in M_r(v)$. We have
$\| v - p(t_n, {\mathtt x}_n) \| \le r$.
From Step~4, we know that
$$
w_n = p(t_n, {\mathtt x}_n) \ \rightharpoonup \
p(t_0, {\mathtt x}_0) = w_0
\ \ \mbox{in} \ \ V,
$$
where $w_0 \in K$ is the unique solution
to problem (\ref{PP4}) corresponding to $(t_0, {\mathtt x}_0)$.
We use the weak lower semicontinuity of the norm to get
$$
\| v - w_0\| \le \liminf \| v - w_n\| =
\liminf \| v - p(t_n, {\mathtt x}_n) \| \le r.
$$
Hence $\| v - p(t_0, {\mathtt x}_0) \| \le r$, so
$(t_0, {\mathtt x}_0) \in M_r(v)$ which proves the closedness
of $M_r(v)$ for all $r> 0$ and $v \in V$. Hence the measurability of (\ref{MEAS1}) follows.

It is known,
see~\cite[Definition~2.5.25, Remark~2.5.26, p. 190]{DMP1}, that for
a measurable space $([0, T], \Sigma)$,
a separable metric space ${\mathbb X}$
and a metric space $V$,
every $\Sigma \otimes {\cal B}({\mathbb X})$ measurable
function is superpositionally measurable.
Therefore, we deduce that the map $p$ is also
superpositionally measurable, that is,
for every $\Sigma$-measurable function
$z \colon [0,T] \to {\mathbb X}$, the function
$$
[0, T] \ni t \ \mapsto \ p(t, z(t)) \in K \subset V
$$
is $\Sigma$-measurable.
Recalling the estimate (\ref{EST33}) in Step~2,
we infer that for any function
$({\bar{\lambda}}, {\bar{\xi}}, {\bar{\eta}}, {\bar{\zeta}})
\in L^2(0, T; E \times X \times Y \times Z)$,
the function
$$
[0, T] \ni t \ \mapsto \ p(t, {\bar\lambda}(t),
{\bar\xi}(t), {\bar\eta}(t), {\bar\zeta}(t)) \in K \subset V
$$
belongs to $L^2(0, T; V)$.
Equivalently, we conclude that for any
$({\bar{\lambda}}, {\bar{\xi}}, {\bar{\eta}}, {\bar{\zeta}})
\in L^2(0, T; E \times X \times Y \times Z)$, then
there exists a unique $w \in L^2(0, T; V)$ such that
$w(t) \in K$ for a.e. $t \in (0, T)$ and
\begin{eqnarray}\label{Limit77}
&&\hspace{-1.0cm}
\langle A (t, {\bar\lambda}(t), w(t)) - f(t,{\bar\xi}(t)),
v - w(t) \rangle
+ \, \varphi (t, {\bar\eta}(t), w(t), v)
- \varphi(t, {\bar\eta}(t), w(t), w(t))
\nonumber \\[1mm]
&&\hspace{-1.0cm}
\ \ \ \quad + \,
j^0 (t, {\bar\zeta}(t), Mw(t); Mv - Mw(t)) \ge 0
\ \ \mbox{\rm for all}
\ \, v \in K, \ \mbox{a.e.} \ t \in (0, T).
\label{MEAS3}
\end{eqnarray}

\noindent
{\bf Step 6}.
(Application of a fixed point argument).

\smallskip

Let $(\lambda_i, \xi_i, \eta_i, \zeta_i) \in
L^2(0,T;E \times X \times Y \times Z)$, $i=1$, $2$
and
$w_1 = w_{\lambda_1\xi_1\eta_1\zeta_1}$,
$w_2 = w_{\lambda_2\xi_2\eta_2\zeta_2} \in L^2(0, T;V)$
with $w_1(t)$, $w_2(t) \in K$ for a.e. $t \in (0,T)$,
be the unique solutions to problem (\ref{Limit77})
corresponding to
$(\lambda_1,\xi_1, \eta_1, \zeta_1)$ and
$(\lambda_2, \xi_2, \eta_2, \zeta_2)$, respectively.
We claim the following estimate holds
\begin{eqnarray}\label{LL888}
&&
\| w_1 - w_2 \|_{L^2(0, t; V)} \le c \,
\big(
\| \lambda_1-\lambda_2\|_{L^2(0,t;E)}
+ \| \xi_1 - \xi_2 \|_{L^2(0, t; X)}
\nonumber
\\[1mm]
&&
\qquad \qquad
+ \, \| \eta_1 - \eta_2 \|_{L^2(0, t; Y)} +
\| \zeta_1 - \zeta_2 \|_{L^2(0, t; Z)}  \big)
\end{eqnarray}
for all $t \in [0,T]$, where $c >0$ is a constant.
From problem~(\ref{Limit77}) it follows that
\begin{eqnarray*}
&&\hspace{-0.6cm}
\displaystyle
\langle A(t, \lambda_1(t), w_1(t))
- f(t,\xi_1 (t)), w_2(t) - w_1(t) \rangle
+ \varphi (t, \eta_1(t), w_1(t), w_2(t))
\\[1mm]
&&\hspace{-0.6cm}
\quad
- \, \varphi(t, \eta_1(t), w_1(t), w_1(t))	
+ j^0 (t, \zeta_1(t), Mw_1(t); Mw_2(t) - Mw_1(t))
\ge 0
\end{eqnarray*}
for a.e. $t \in (0, T)$
and
\begin{eqnarray*}
&&\hspace{-0.6cm}
\displaystyle
\langle A(t, \lambda_2(t), w_2(t))
- f(t,\xi_2 (t)), w_1(t) - w_2(t) \rangle
+ \varphi (t, \eta_2(t), w_2(t), w_1(t))
\\[1mm]
&&\hspace{-0.6cm}
\quad
- \, \varphi(t, \eta_2(t), w_2(t), w_2(t))	
+  j^0 (t, \zeta_2(t), Mw_2(t); Mw_1(t) - Mw_2(t))
\ge 0	
\end{eqnarray*}
for a.e.\ $t \in (0, T)$.
Next, we add the last two inequalities to get
\begin{eqnarray*}
&&\hspace{-0.9cm}
\langle A(t, \lambda_1(t), w_1(t))
- A(t, \lambda_2(t), w_2(t)), w_1(t) - w_2(t)\rangle
\\ [1mm]
&&\hspace{-0.5cm}
\le
\langle f(t, \xi_1(t)) - f(t, \xi_2(t)),
w_1(t) - w_2(t) \rangle
+ \varphi(t, \eta_1(t), w_1(t), w_2(t))
\\[1mm]
&&\hspace{-0.3cm}
- \, \varphi (t, \eta_1(t), w_1(t), w_1 (t))
+ \, \varphi(t, \eta_2(t), w_2(t), w_1(t))
- \varphi (t, \eta_2(t), w_2(t), w_2 (t)) \\[1mm]
	&&\hspace{-0.1cm}
	+ \, j^0(t, \zeta_1(t), Mw_1(t); Mw_2(t) - Mw_1(t))
	+ j^0(t, \zeta_2(t), Mw_2(t); Mw_1(t) - Mw_2(t))
\end{eqnarray*}
for a.e. $t \in (0,T)$.
By an argument exploiting $H(A)$(d), $H(f)$(b), $H(\varphi)$(b) and $H(j)$(d),
we obtain
\begin{eqnarray*}
&&\hspace{-0.5cm}
m_A \, \| w_1(t) - w_2(t) \|^2
\\[2mm]
&&\hspace{-0.3cm}
\le
{\bar{m}}_A \| \lambda_1(t)-\lambda_2(t)\|_E \,
\| w_1(t)-w_2(t)\|
+ L_f \| \xi_1(t)-\xi_2(t)\|_{X} \, \|w_1(t)-w_2(t)\|
\\[2mm]
&&\hspace{-0.1cm}
+ \, \alpha_{\varphi} \, \| w_1(t) - w_2(t) \|^2
+ \beta_{\varphi} \| \eta_1(t)-\eta_2 (t)\|_Y
\, \| w_1(t) - w_2(t) \|
\\[2mm]
&&\
+ \, m_j \, \| Mw_1(t) - Mw_2(t) \|_X^2
+ \, m_1 \| \zeta_1(t)-\zeta_2(t)\|_Z \| M \|
\| w_1 (t) - w_2(t) \|
\end{eqnarray*}
for all $t \in [0,T]$, and subsequently
\begin{eqnarray*}
&&\hspace{-0.4cm}
( m_A - m_j \|M_1\|^2 -\alpha_{\varphi}) \,
\| w_1 (t)- w_2 (t)\|^2
\le \Big(
{\bar{m}}_A \| \lambda_1(t)-\lambda_2(t)\|_E
+ L_f \| \xi_1(t)-\xi_2(t)\|_{X} \\
&&\qquad\quad
+ \, \beta_{\varphi} \, \| \eta_1(t)-\eta_2 (t)\|_Y
+ \, m_1 \| M \| \| \zeta_1(t)-\zeta_2(t)\|_Z \Big)
\| w_1 (t) - w_2(t) \|
\end{eqnarray*}
for all $t \in [0,T]$.
We integrate the latter on $(0, t)$ for all $t \in [0,T]$, and
use the H\"older inequality to get
\begin{eqnarray*}
&&\hspace{-0.4cm}
( m_A - m_j \|M_1\|^2 -\alpha_{\varphi}) \,
\| w_1 - w_2 \|^2_{L^2(0,t;V)}
\le c\, \Big(
\| \lambda_1-\lambda_2\|_{L^2(0,t;E)}
+ \| \xi_1-\xi_2\|_{L^2(0,t;X)} \\
&&\qquad\quad
+ \, \| \eta_1-\eta_2 \|_{L^2(0,t;Y)}
+ \, \| \zeta_1 -\zeta_2 \|_{L^2(0,t;Z)} \Big)
\| w_1  - w_2 \|_{L^2(0,t;V)}
\end{eqnarray*}
for all $t \in [0,T]$ with $c > 0$.
By $(H_0)$ we deduce that the inequality
(\ref{LL888}) holds.

\medskip

We are now in a position to apply a fixed point argument.
We define the operator
$\Lambda \colon L^2(0,T; E \times X \times Y \times Z)
\to L^2(0,T;E \times X \times Y \times Z)$ by
$$
\Lambda(\lambda, \xi, \eta, \zeta) =
(R_1 w_{\lambda\xi\eta\zeta},
{R}_2 w_{\lambda\xi\eta\zeta},
{R_3} w_{\lambda\xi\eta\zeta},
{R_4} w_{\lambda\xi\eta\zeta})
$$
for all
$(\lambda, \xi, \eta, \zeta) \in L^2(0,T;E \times X \times Y \times Z)$,
where $w_{\lambda\xi\eta\zeta} \in L^2(0, T; V)$
denotes the unique solution
to problem~(\ref{Limit77}) corresponding to
$(\lambda, \xi, \eta, \zeta)$.
We exploit hypothesis $H(R)$,
inequality (\ref{LL888}), the H\"older inequality,
and similarly as in~\cite[Theorem~3, Step~4]{MBZ},
we find a constant $c > 0$ such that
\begin{eqnarray}
&&\| \Lambda(\lambda_1,\xi_1, \eta_1, \zeta_1)(t) -
\Lambda(\lambda_2, \xi_2, \eta_2, \zeta_2)(t)
\|^2_{E\times X \times Y \times Z}
\nonumber \\
&&\qquad\qquad
\le c \, \int_{0}^{t} \,
\| (\lambda_1, \xi_1, \eta_1, \zeta_1)(s) -
(\lambda_2, \xi_2, \eta_2, \zeta_2)(s)
\|^2_{E \times X \times Y \times Z} \, ds
\label{pstaly}
\end{eqnarray}
for all $t \in [0,T]$.
From~\cite[Corollary~27]{SM2},
we infer that there exists a unique fixed point
$(\lambda^*, \xi^*,\eta^*,\zeta^*)$ of $\Lambda$,
that is
$$
(\lambda^*, \xi^*, \eta^*,\zeta^*) \in
L^2(0,T;E \times X \times Y \times Z)
\ \ {\rm with}\ \
\Lambda(\lambda^*, \xi^*, \eta^*, \zeta^*)
= (\lambda^*, \xi^*, \eta^*, \zeta^*).
$$

Given $(\lambda^*, \xi^*, \eta^*, \zeta^*) \in
L^2(0,T; E \times X \times Y \times Z)$
the unique fixed point of the operator $\Lambda$,
we define $w_{\lambda^* \xi^* \eta^* \zeta^*}
\in L^2(0, T; V)$
to be the unique solution to problem~(\ref{Limit77}) corresponding to $(\lambda^*, \xi^*,\eta^*, \zeta^*)$.
By the definition of the operator $\Lambda$, we have
$$
\lambda^* = R_1(w_{\lambda^* \xi^* \eta^* \zeta^*}),
\ \
\xi^* = {R}_2(w_{\lambda^* \xi^* \eta^* \zeta^*}),
\ \
\eta^* = {R_3}(w_{\lambda^* \xi^* \eta^* \zeta^*})
\ \ \mbox{and} \ \
\zeta^* = {R_4}(w_{\lambda^* \xi^* \eta^* \zeta^*}) .
$$
Finally, we use these relations in problem~(\ref{Limit77}), and conclude that
$w_{\lambda^* \xi^* \eta^* \zeta^*}\in L^2(0, T; V)$
is the unique solution to Problem~\ref{Problem1q}.
This completes the proof of the theorem.
\hfill$\Box$

\section{Application to a quasistatic contact problem}
\label{Application}

\noindent
Several quasistatic contact problems in solid mechanics lead
to a variational-hemivariatio\-nal inequality involving history-dependent operators of the form $(*)$
in which the unknown is the velocity field.
We illustrate the applicability of results of Section~\ref{s2}
to a unilateral viscoelastic frictional contact problem to which Theorem~\ref{theorem1q} can be applied.

We begin with the physical setting of the model.
The reference configuration of a viscoelastic
body is a bounded domain $\Omega$ in $\real^d$,
$d=2$, $3$ with Lipschitz boundary such that
$\partial \Omega = \Gamma_1 \cup \Gamma_2 \cup \Gamma_3 \cup \Gamma_4$
with mutually disjoint and measurable parts $\overline{\Gamma}_{1}$, $\overline{\Gamma}_{2}$, $\overline{\Gamma}_{3}$ and
$\overline{\Gamma}_{4}$ such that
${\rm meas}\,(\Gamma_{1}) > 0$.
In what follows, $\bnu$ denotes the outward unit normal at
the boundary and $\mathbb{S}^{d}$ stands for the
space of $d\times d$ symmetric matrices.
The classical contact model under consideration is represented the following boundary value problem.
\begin{Problem}\label{CONTACTq}
	Find a displacement field
	$\bu \colon \Omega\times (0, T) \to\mathbb{R}^d$ and a stress field
	$\bsigma \colon \Omega\times (0, T) \rightarrow \mathbb{S}^d$ such that for all $t\in (0,T)$,
	\begin{align}
	\label{equation1}
	\bsigma(t)
	&={\mathscr A}\bvarepsilon ({\bu}'(t))
	+{\mathscr B}\bvarepsilon({\bu}(t))
	+\int_0^t{\mathscr C}(t-s)\bvarepsilon({\bu}'(s))\,ds \quad
	&{\rm in}\
	&\Omega,\\[1mm]
	\label{equation2}
	{\rm Div} \, \bsigma (t)&+\fb_0(t) = 0
	\quad&{\rm in}\ &\Omega,\\[1mm]
	\label{equation3} \bu(t)&=\bzero &{\rm on}\ &\Gamma_1,\\[1mm]
	\label{equation4} \bsigma(t)\bnu&=\fb_N(t)\quad&{\rm on}\ &\Gamma_2,\\[1mm]
	\nonumber
	u'_{\nu}(t)&\le g, \
	\sigma_{\nu}(t)+\eta(t) \le 0, \ (u'_{\nu}(t)-g)(\sigma_{\nu}(t)+\eta(t))=0,
	\\[1mm]
	\label{equation5}
	\eta(t)&\in k(u_{\nu}(t))\,
	\partial j_{\nu}(u'_{\nu}(t))
	\quad&{\rm	on}\
	&\Gamma_3,\\[1mm]
	\|\bsigma_\tau(t)\|&\le F_b\Big(t,\int_0^t\|\bu_\tau(s)\|\,ds\Big),\quad &\ &\nonumber\\[1mm]
	-\bsigma_\tau(t)
	&=F_b\Big(t,\int_0^t\|\bu_\tau(s)\|\,ds\Big)
	\frac{{\bu}_\tau'(t)}{\|{\bu}_\tau'(t)\|}\ \ \  {\rm  if}\
	{\bu}_\tau'(t)\ne\bzero \quad&{\rm on}\ &\Gamma_3,\label{equation6}
	\end{align}
\begin{align}\label{E429}
\quad -\sigma_\nu(t)&=p_\nu({u}_\nu(t)),
&{\rm on}\ &\Gamma_4,
& \\[2mm]
&\hspace{-1.3cm}
\begin{cases}
\|\bsigma_\tau(t)\|&\le \mu (\|\bu'_\tau (t)\|)
|\sigma_\nu (t)|, \\[2mm]
\displaystyle
- \bsigma_\tau (t)&
\displaystyle= \mu (\|\bu'_\tau (t)\|) |\sigma_\nu(t)|
\frac{{\bu}'_\tau(t)}{\|{\bu}'_\tau(t)\|}
\  {\rm  if}\
{\bu}'_\tau(t)\ne\bzero \ 
\end{cases}
&{\rm on}\ &\Gamma_4,
\label{E430}
\end{align}
	and
	\begin{equation}\label{equation7}
	\bu(0)=\bu_0 \ \ {\rm in}\quad \Omega.
	\end{equation}
\end{Problem}

\smallskip

\noindent
In this problem,
the normal and tangential components on the boundary of
a vector $\bv$ are defined by
$v_{\nu}=\bv\cdot \bnu$ and $\bv_{\tau}=\bv-v_{\nu}\bnu$,
respectively.
Given a tensor $\bsigma$, the symbols
$\sigma_{\nu}$ and $\bsigma_{\tau}$ stand for
its normal and tangential components on the boundary,
that is,
$\sigma_{\nu}=(\bsigma\bnu)\cdot\bnu$
and
$\bsigma_{\tau}=\bsigma\bnu-\sigma_{\nu}\bnu$.
The linearized strain tensor defined by
$$
\bvarepsilon(\bu) = (\varepsilon_{ij}(\bu)), \ \ \
\varepsilon_{ij} (\bu) = \frac{1}{2} (u_{i,j} + u_{j,i})
\ \ \mbox{in} \ \ \Omega.
$$

We now present a description of conditions in Problem~\ref{CONTACTq}
together with the hypotheses on the data.
Equation (\ref{equation1}) is the constitutive
law for viscoelastic materials with long memory in which $\mathscr{A}$ denotes the viscosity operator,
$\mathscr{B}$ represents the elasticity operator,
$\mathscr{C}$ is the relaxation tensor, and
which satisfy the following hypotheses.

\smallskip

\noindent
$\underline{H({\mathscr{A}})}:$ \quad
${\mathscr{A}} \colon \Omega\times [0, T] \times  \mathbb{S}^{d}\rightarrow \mathbb{S}^{d}$
is the viscosity operator such that

\smallskip

\lista{
	\item[(1)]
	${\mathscr{A}}(\cdot,t,\bvarepsilon)$
	is continuous for all $t \in [0, T]$,
	$\bvarepsilon \in \mathbb{S}^{d}$.
	\smallskip
	\item[(2)]
	${\mathscr{A}}(\bx,\cdot,\cdot)$
	is continuous for a.e. $\bx \in \Omega$.
	\smallskip
	\item[(3)]
	$\| {\mathscr{A}} (\bx, t, \bvarepsilon) \|
	\le {\widetilde{a}}_0(\bx, t) + {\widetilde{a}}_1
	\| \bvarepsilon \|$
	for all $\bvarepsilon \in \mathbb{S}^{d}$,
	$t \in [0,T]$, a.e. $\bx \in \Omega$ \\
	with ${\widetilde{a}}_0 \in C(0, T; L^2(\Omega))_+$,
	${\widetilde a}_1 > 0$.
	\smallskip
	\item[(4)]
	$\mbox{there exists} \ m_{\mathscr{A}}>0 \ \mbox{such that} \ \
	(\mathscr{A}(\bx,t,\bvarepsilon_{1})-
	\mathscr{A}(\bx,t,\bvarepsilon_{2}))\cdot (\bvarepsilon_{1}-\bvarepsilon_{2})\ge m_{\mathscr{A}}
	\|\bvarepsilon_{1}-\bvarepsilon_{2}\|^{2}$ \\
	$\ \ \qquad \mbox{for all} \  \bvarepsilon_{1}, \bvarepsilon_{2}\in \mathbb{S}^{d}, \ \mbox{all} \ t \in [0, T], \, \mbox{a.e.} \ \bx\in \Omega$.
}

\smallskip

\noindent
$\underline{H({\mathscr B})}:$ \quad
$\mathscr {B}\colon \Omega\times  \mathbb{S}^{d}\rightarrow \mathbb{S}^{d}$
is the elasticity operator such that

\smallskip

\lista{
	\item[(1)]
	$\mathscr{B}(\cdot,\bvarepsilon)$
	is measurable on $\Omega$ for all
	$\bvarepsilon \in \mathbb{S}^{d}$.
	\smallskip
	\item[(2)]
	there exists $L_{\mathscr{B}}>0$ such that
	$\|\mathscr{B}(\bx,\bvarepsilon_{1})-
	\mathscr{B}(\bx,\bvarepsilon_{2})\|\leq L_{\mathscr{B}}
	\|\bvarepsilon_{1}-\bvarepsilon_{2}\|$
	for all \\ $\bvarepsilon_{1}$, $\bvarepsilon_{2}\in \mathbb{S}^{d}$, a.e. $\bx \in \Omega$.
	\smallskip
	\item[(3)]
	${\mathscr B}(\bx,  {\bf{0}}) = {\bf{0}}$
	for a.e. $\bx \in \Omega$.
}

\smallskip

\noindent
$\underline{H({\mathscr C})}:$ \quad
$\mathscr{C}\in C(0,T; Q_\infty)$ is the relaxation  operator.

\smallskip

In $H({\mathscr C})$, the space $Q_\infty$
of fourth order tensors is defined by
\begin{equation*}
Q_\infty=\{\, \bsigma=(\sigma_{ijkl}) \mid \sigma_{ijkl}=\sigma_{jikl}
=\sigma_{klij}\in L^{\infty}(\Omega),
1\le i,j,k,l\le d \, \},
\end{equation*}
which a Banach space with the norm
\begin{equation*}
\|\bsigma\|_{Q_\infty}
=\sum_{1\le i,j,k,l\le d}
\|\sigma_{ijkl}\|_{L^{\infty}(\Omega)}
\ \ \mbox{for} \ \ \bsigma \in Q_\infty.
\end{equation*}

Equation (\ref{equation2}) represents
the equilibrium equation in which
${\rm Div} \bsigma = (\sigma_{ij,j})$ and
$\fb_{0}$ denotes the density of the body forces.
The boundary condition (\ref{equation3}) states that the displacement vanishes, which means that the body is fixed along  $\Gamma_1$. Relation (\ref{equation4}) is the traction boundary condition with surface tractions of density $\fb_N$ acting on the part $\Gamma_2$ of the boundary.
Condition (\ref{equation5}) is the Signorini  unilateral contact boundary condition for the normal velocity
in which $\partial j_\nu$ stands for the Clarke
subgradient of a prescribed function $j_\nu$ and $g > 0$ is a constant.
Condition $\eta(t) \in k(u_{\nu}(t))\,
\partial j_{\nu}(u'_{\nu}(t))$ on $\Gamma_3$ represents a generalization of the normal damped response
condition where $k$ is a given damper
coefficient depending on the normal displacement.
The assumptions on $k$ and $j_\nu$ read as follows.

\smallskip

\noindent
$\underline{H({k})}:$ \quad
$k \colon \Gamma_3 \times [0, T] \times \real \times \real \to \real$
is a damper coefficient such that

\smallskip

\lista{
	\item[(1)]
	$k(\cdot, t, r)$ is continuous for all
	$t \in [0, T]$, $r \in \real$.
	\smallskip
	\item[(2)]
	$k(\bx, \cdot, r)$ is continuous for all
	$r \in \real$, a.e. $\bx \in \Gamma_{3}$.
	\smallskip
	\item[(3)]
	there exist $0<k_1\le k^*$ such that
	$0 < k_1 \le k(\bx,t, r)\le k^*$ for all
	$t \in [0, T]$, $r \in \real$,
	a.e. $\bx \in \Gamma_{3}$.
	\smallskip
	\item[(4)]
	there exists $L_{k}>0$ such that  $|k(\bx,t,r_{1})-k(\bx,t,r_{2})|\le L_{k}|r_{1}-r_{2}|$
	for all $t \in [0, T]$, $r_{1}$, $r_{2}\in \real$,
	a.e. $\bx \in \Gamma_{3}$.
}

\smallskip

\noindent
$\underline{H({j_\nu})}:$ \quad
$j_\nu \colon \real \to \real$
is a function such that

\smallskip

\lista{
	\item[(1)]
	$j_{\nu}$ is locally Lipschitz, and
	either $j_\nu$ or $-j_\nu$ is Clarke regular.
	\smallskip
	\item[(2)]
	there exists $\overline{c}_{0} >0$	such that
	$|\partial j_{\nu}(r)|\le \overline{c}_{0}$
	for all $r \in \real$.
	\smallskip
	\item[(3)]
	there exists $m_{j_{\nu}} \ge 0$ such that
	$j_{\nu}^0(r_{1}; r_{2}-r_{1}) + j_{\nu}^0(r_{2};r_{1}-r_{2})
	\le m_{j_{\nu}}|r_{1}-r_{2}|^{2}$
	for all $r_{1}$, $r_{2}\in \mathbb{R}$.
}

\smallskip

Condition (\ref{equation6}) represents
a version of the Coulomb law of dry friction
in which the function $F_b$ denotes the friction bound.
The latter may depend on the accumulated (total) slip
represented by the quantity
$
\int_0^t \| \bu_{\tau}(\bx,s) \| \, ds
$
at the point $\bx \in \Gamma_3$
in the time interval $[0,t]$.
Conditions (\ref{E429}) and (\ref{E430}) represent the contact condition with normal compliance and the Coulomb law of dry friction, respectively.
Here, $p$ is the normal compliance functions and
$\mu$ denotes the coefficient of friction
depending on the slip-rate $\bu'_\tau$
being the tangential part of the velocity field.
The hypotheses on the data in
(\ref{equation6})--(\ref{E430}) are the following.

\smallskip

\noindent
$\underline{H({F_b})}:$ \quad
$F_b \colon \Gamma_3 \times [0, T] \times \real \to
[0, \infty)$ is the friction bound such that

\smallskip

\lista{
	\item[(1)]
	$F_b(\cdot, t, r)$ is measurable
	for all $t \in [0, T]$, $r \in \real$.
	\smallskip
	\item[(2)]
	$F_b(\bx, \cdot, r)$ is continuous
	for all $r \in \real$, a.e. $\bx \in \Gamma_3$.
	\item[(3)]
	there exists $L_{F_b}>0$ such that
	$|F_b(\bx,t,r_{1}) - F_b(\bx,t,r_{2})|\le L_{F_b}|r_{1}-r_{2}|$ for all $t \in [0, T]$,
	$r_{1}$, $r_{2}\in \mathbb{R}$, a.e. $\bx \in \Gamma_{C}$.
	\smallskip
	\item[(4)]
	$\bx \mapsto F_b(\bx,t, 0)$ belongs to $L^2(\Gamma_3)$
	for all $t \in [0, T]$.
}

\smallskip

\noindent
$\underline{H({p})}:$ \quad
$p \colon \Gamma_4 \times [0, T] \times \real \to
[0, \infty)$ is the normal compliance function
such that

\smallskip

\lista{
	\item[(1)]
	$p(\cdot, t, r)$ is continuous for all
	$t \in [0, T]$, $r \in \real$.
	\smallskip
	\item[(2)]
	$p(\bx, \cdot, r)$ is continuous for all
	$r \in \real$, a.e. $\bx \in \Gamma_4$.
	\smallskip
	\item[(3)]
	there exists $p^* > 0$ such that $p(\bx, t, r) \le p^*$
	for all $t \in [0, T]$,	$r \in \real$,
	a.e. $\bx \in \Gamma_{4}$.
	\smallskip
	\item[(4)]
	$p(\bx, t, r) = 0$ for all $t \in [0, T]$,	$r < 0$,
	a.e. $\bx \in \Gamma_{4}$.
	\smallskip
	\item[(5)]
	there exists $L_p>0$ such that  $|p(\bx,t,r_{1})-k(\bx,t,r_{2})|\le L_{p}|r_{1}-r_{2}|$
	for all $t \in [0, T]$, $r_{1}$, $r_{2}\in \real$,
	a.e. $\bx \in \Gamma_{4}$.
}

\smallskip

\noindent
$\underline{H({\mu})}:$ \quad
$\mu \colon \Gamma_4 \times [0, T] \times [0, \infty)
\to [0, \infty)$ is the coefficient of friction
such that

\smallskip

\lista{
	\item[(1)]
	$\mu(\cdot, t, r)$ is continuous for all
	$t \in [0, T]$, $r \in \real$.
	\smallskip
	\item[(2)]
	$\mu(\bx, \cdot, r)$ is continuous for all
	$r \in \real$, a.e. $\bx \in \Gamma_4$.
	\smallskip
	\item[(3)]
	there exists $\mu^* > 0$ such that
	$\mu(\bx, t, r) \le p^*$
	for all $t \in [0, T]$,	$r \in \real$,
	a.e. $\bx \in \Gamma_{4}$.
	\smallskip
	\item[(4)]
	there exists $L_\mu>0$ such that  $|\mu(\bx,t,r_{1})-\mu(\bx,t,r_{2})|\le L_{\mu}|r_{1}-r_{2}|$
	for all $t \in [0, T]$,
	$r_{1}$, $r_{2}\in [0, \infty)$,
	a.e. $\bx \in \Gamma_{4}$.
}

The normal compliance contact condition was first introduced in~\cite{Oden} and used in many publications,
see, e.g. \cite{HS,MOSBook,SST,SMatei,SM2}.
%
Details on mechanical interpretation of conditions (\ref{equation5}) and (\ref{equation6}) can be found in~\cite{hms,MOSBook} and references therein while for conditions (\ref{E429}) and (\ref{E430}) we refer to~\cite{SXiao2015}. Finally,
the densities of the body forces and the surface tractions,
and the initial displacement in (\ref{equation7})
are supposed to have the following regularity.

\medskip

\noindent
$\underline{(H_1)}:$ \quad
$\fb_0 \in C(0,T;L^2(\Omega;\real^d))$,
$\fb_N \in C(0,T;L^2(\Gamma_2;\real^d))$,
$\bu_0 \in V$.

\medskip

Let
$\mathcal{H}=L^{2}(\Omega;\mathbb{S}^{d})$
be the Hilbert space with the standard inner product
\begin{equation*}
\langle \bsigma, \bvarepsilon \rangle_{\mathcal{H}}=
\int_{\Omega}\sigma_{ij}(\bx)\varepsilon_{ij}(\bx)\, dx
\ \ \mbox{for all} \ \ \bsigma, \bvarepsilon \in {\cal H}.
\end{equation*}
Consider also the space
$V=\{\, \bv\in H^{1}(\Omega;\mathbb{R}^{d})\mid \bv={\bf{0}}~\textrm{on}~\Gamma_{1}\, \}$ equipped
with the inner product
$(\bu,\bv) = (\bvarepsilon(\bu),
\bvarepsilon(\bv))_{\mathcal{H}}$ and the
corresponding norm
$\|\bv\| = \| \bvarepsilon(\bv)\|_{\mathcal{H}}$
for all $\bu$, $\bv\in V$.
It is well known that
$V$ is a Hilbert space and there is a linear trace
operator
$\gamma\colon V\to L^{2}(\Gamma;\mathbb{R}^{d})$
such that
$\|\bv \|_{L^{2}(\Gamma;\mathbb{R}^{d})}
\le \|\gamma\|\|\bv \|$ for all $\bv\in V$,
where $\|\gamma\|$ denotes the norm of the trace
operator in ${\cal L}(V,L^2(\Gamma;\real^d))$.
%
%
Further, we introduce the set of admissible velocity fields
$U$ defined by
\begin{equation*}
U=\{\, \bv \in V \mid
v_{\nu} \leq g \ \ \mbox{\rm on} \ \ \Gamma_{3}\, \},
\end{equation*}
and a linear bounded functional on $V$
\begin{equation}\label{fXXX}
\bv \mapsto \langle \fb (t),\bv\rangle
=
\langle \fb_{0}(t), \bv\rangle_{L^{2}(\Omega;\mathbb{R}^{d})}
+\langle \fb_{N}(t),\bv \rangle_{L^{2}(\Gamma_{2};\mathbb{R}^{d})}
\end{equation}
for all $\bv \in V$. From $(H_1)$ this functional has the regularity $\fb \in C(0, T; V^*)$.
In what follows we often do not indicate explicitly the dependence of functions and operators on the variable $\bx$.

We now shortly sketch the procedure to obtain
the weak formulation of Problem~\ref{CONTACTq}.
Let $(\bu, \bsigma)$ be a smooth solution to this problem which means that the data are smooth functions such that all the derivatives and all the conditions are satisfied in the usual sense at each point.
Let $\bv\in U$ and $t\in [0,T]$.
We multiply the equilibrium equation (\ref{equation2})
by $\bv-\bu'(t)$,
use  the integration by parts formula,
and apply the boundary conditions (\ref{equation3})
and (\ref{equation4})
to deduce
\begin{eqnarray*}
&&\hspace{-0.5cm}
\int_\Omega \bsigma(t)\cdot
\big(
\bvarepsilon(\bv) -\bvarepsilon(\bu'(t))
\big) \, dx =\int_\Omega \fb_0(t)\cdot(\bv-\bu'(t))dx
\\ [1mm]
&&\hspace{-0.5cm} \qquad \ \
+\int_{\Gamma_2} \fb_N(t)\cdot(\bv-\bu'(t))\, d\Gamma
+\int_{\Gamma_3\cup\Gamma_4} \bsigma(t)\bnu \cdot(\bv-\bu'(t))\, d\Gamma.
\end{eqnarray*}
The unilateral contact condition (\ref{equation5})
and the definition of the Clarke subgradient imply that
\begin{eqnarray}
&&\hspace{-2.0cm}
\sigma_{\nu}(t) (v_\nu - u_\nu'(t))
= (\sigma_{\nu}(t) + \eta(t))(v_\nu - g)
- (\sigma_{\nu}(t) + \eta(t))(u_\nu'(t) - g)
\nonumber \\[2mm]
&&\qquad \
- \, \eta(t) (v_\nu - u_\nu'(t))
\ge
- k(u_\nu(t)) j_\nu^0(u_\nu'(t); v_\nu - u_\nu'(t))
\ \ \mbox{on} \ \ \Gamma_3 ,
\label{normal}
\end{eqnarray}
while the friction law (\ref{equation6})
can be equivalently described by
\begin{equation}\label{tangent}
\bsigma_{\tau}(t) \cdot (\bv_\tau - \bu_\tau'(t))
\ge
- F_{b}\Big(t, \int_0^t\|\bu_{\tau}(s)\|\,ds\Big)
\big(
\|\bv_\tau \|- \|\bu_\tau'(t)\| \big)
\ \ \mbox{on} \ \ \Gamma_3.
\end{equation}
We combine (\ref{normal}), (\ref{tangent}) and invoke
the decomposition formula in~\cite[(6.33)]{MOSBook}
to see that
\begin{eqnarray*}
	&&F_{b}\Big(t, \int_0^t\|\bu_{\tau}(s)\|\,ds\Big)
	\big(
	\|\bv_{\tau}\|-\|\bu'_{\tau}(t)\| \big)
	+k(u_{\nu}(t)) \, j_{\nu}^0(u'_{\nu}(t);
	v_{\nu}-u'_{\nu}(t))\\ [1mm]
	&&\qquad\qquad
	+\, \bsigma(t)\bnu \cdot(\bv-\bu'(t))\ge 0
	\ \ \mbox{on} \ \ \Gamma_3.
\end{eqnarray*}
From (\ref{E429}) and (\ref{E430}),
we obtain the boundary integrals on the part $\Gamma_4$.
Hence and by the definition (\ref{fXXX}) we get
\begin{eqnarray*}
&&
\langle \bsigma(t), \bvarepsilon(\bv)-\bvarepsilon(\bu'(t))
\rangle_{\mathcal{H}}
+ \int_{\Gamma_{3}} F_{b}\Big(t, \int_0^t\|\bu_{\tau}(s)\| \,ds\Big)
(\|\bv_{\tau}\| - \|\bu'_{\tau}(t)\|)\, d\Gamma
\\ [1mm]
&& \qquad +\int_{\Gamma_{3}} k(u_{\nu}(t))\,j_{\nu}^0(u'_{\nu}(t);
v_{\nu}-u'_{\nu}(t))\, d\Gamma
+
\int_{\Gamma_4}
p(u_\nu(t)) \, (v_\nu - u'_\nu(t)) \, d\Gamma
\\[1mm]
&&\qquad\qquad
+ \int_{\Gamma_4}
\mu (\|\bu'_\tau (t)\|) \, p(u_\nu(t)) \,
(\|\bv_\tau \|- \|\bu'_\tau(t)\| ) \, d\Gamma
\ge
\langle \fb(t), \bv - \bu \rangle.
\end{eqnarray*}
Finally, we use the constitutive relation (\ref{equation1}),
we arrive to the following variational formulation of Problem~\ref{CONTACTq}.
\begin{Problem}\label{CONTACT1q}
Find $\bu \colon (0,T)\to V$ such that $\bu(0)=\bu_0$
and
\begin{eqnarray*}
&&\hspace{-0.5cm}
\langle \mathscr{A}(\bvarepsilon(\bu'(t)))
+ \mathscr{B}(\bvarepsilon(\bu(t)))
+ \int_0^t \mathscr{C}(t-s)\bvarepsilon(\bu'(s))\, ds, \bvarepsilon(\bv)-\bvarepsilon(\bu'(t))
\rangle_{\mathcal{H}} \\
&&
\hspace{-0.4cm}
\ \ + \int_{\Gamma_{3}}
F_{b}\Big(t, \int_0^t\|\bu_{\tau}(s)\| \, ds \Big)
(\|\bv_{\tau}\| - \|\bu'_{\tau}(t)\|) \, d\Gamma
+\int_{\Gamma_{3}}
k(u_{\nu}(t)) \,
j_{\nu}^0(u'_{\nu}(t); v_{\nu}-u'_{\nu}(t))\, d\Gamma
\\[2mm]
&&
+
\int_{\Gamma_4}
p(u_\nu(t)) \, (v_\nu - u'_\nu(t)) \, d\Gamma
+ \int_{\Gamma_4}
\mu (\|\bu'_\tau (t)\|) \, p(u_\nu(t)) \,
(\|\bv_\tau \|- \|\bu'_\tau(t)\| ) \, d\Gamma \\[2mm]
&&\qquad
\ge \langle \fb(t), \bv - \bu \rangle
\ \ \mbox{\rm for all} \ \ \bv \in U,
\ \mbox{\rm all} \ t \in [0,T].
\end{eqnarray*}
\end{Problem}

The following result concerns the unique solvability and regularity of solution to Problem~\ref{CONTACT1q}.
\begin{Theorem}\label{existence3}
Assume hypotheses
$H({\mathscr A})$, $H({\mathscr B})$,
$H({\mathscr C})$, $H(F_b)$, $H(p)$, $H(\mu)$,
$H(k)$, $H(j_\nu)$, $(H_1)$,
and the following smallness condition
	\begin{equation}\label{CONTACT_smallness}
	k^* m_{j_\nu} \| \gamma\|^4 + p^* L_\mu \| \gamma \|^2 <
	m_{\mathscr{A}}.
	\end{equation}
	Then Problem~$\ref{CONTACT1q}$ has a unique solution
	$\bu \in C(0, T; V)$ with $\bu' \in L^2(0, T; V)$ and
	$\bu'(t) \in U$ for a.e. $t \in (0, T)$.
\end{Theorem}

\noindent
{\bf Proof}. \
We will apply Theorem~\ref{theorem1q}
with the following functional framework:
$E = {\cal H}$, $X=L^2(\Gamma_3)$,
$Y = L^2(\Gamma_{3}) \times L^2(\Gamma_{4})$,
$Z=L^2(\Gamma_{3})$ and $K=U$.

Let the operator
$A \colon [0, T]\times E \times V \to V^*$, and functions
$\varphi \colon [0, T] \times Y \times V \times V
\to \mathbb{R}$
and
$j \colon [0, T] \times Z\times X\to \mathbb{R}$
be defined by
\begin{eqnarray*}
&&
\langle A(t, \lambda, \bv), \bz \rangle
= \langle \mathscr{A}(t, \bvarepsilon(\bv)) + \lambda,
\bvarepsilon(\bz) \rangle_{\mathcal{H}}
\ \ \mbox{for} \ \ t \in [0, T], \, \lambda \in E,
\, \bv, \bz \in V,
\\ [1mm]
&&
\varphi(t, \eta, \bw, \bv)
= \int_{\Gamma_{3}} F_{b}(t, \eta_1) \| \bv_\tau \|
\, d\Gamma
+ \int_{\Gamma_4}
\left(
p(t, \eta_2) v_\nu + \mu(t, \| \bw_\tau \|) \,
p(t, \eta_2)  \| \bv_\tau \| \right) \, d\Gamma
\\[1mm]
&&\qquad\qquad
\ \ \mbox{for} \ \ t \in [0, T],
\eta= (\eta_1, \eta_2) \in Y, \bw, \bv \in V,
\\ [1mm]
&&
j(t, \zeta, \bv)= \int_{\Gamma_{3}} k(t,\zeta) \, j_{\nu}(v)\, d\Gamma
\ \ \mbox{for} \ \ t \in [0, T], \, \zeta \in Z,
\bv \in V.
\end{eqnarray*}
We introduce operators
$I \colon L^2(0, T; V) \to L^2(0, T; V)$,
$R_1 \colon L^2(0, T; V) \to L^2(0, T; E)$,
$R_2 = 0$,
$R_3 \colon L^2(0, T; V) \to L^2(0, T; Y)$, and
$R_4 \colon L^2(0, T; V) \to L^2(0,T;Z)$
given by
\begin{eqnarray*}
&&
(I\bw)(t) = \bu_0 + \int_0^t \bw (s)\,ds
\ \ \mbox{for} \ \ \bw \in L^2(0, T; V), \, t \in [0, T],
\\
&&
(R_1\bw)(t) =
\mathscr{B} \bvarepsilon ((I\bw)(t)) +
\int_0^t {\mathscr C}(t-s) \bvarepsilon (\bw(s))
\, ds
\ \ \mbox{for} \ \ \bw \in L^2(0, T; V), \, t \in [0, T],
\\
&&(R_3\bw)(t)=
\Big(
\int_0^t \|\int_0^s\bw_\tau(r)dr+\bu_{0\tau}\| \, ds,
((I\bw)(t))_\nu \Big)
\  \mbox{for} \  \bw \in L^2(0, T; V), \,
t \in [0, T], \\[1mm]
&&
(R_4\bw)(t)= ((I\bw)(t))_\nu
\ \ \mbox{for} \ \ \bw \in L^2(0, T; V),\, t \in [0, T].
\end{eqnarray*}
We also define operator $M \colon V \to X$ by
$M \bv = v_\nu$ for $\bv \in V$.
Let $\bw(t)=\bu'(t)$ for all $t\in [0,T]$.
Then, with the above notation, we consider the inequality problem associated with Problem~\ref{CONTACT1q}.
\begin{Problem}\label{CONTACT2q}
Find $\bw\in L^2(0, T; V)$ such that
$\bw(t) \in U$ for a.e. $t \in (0, T)$
and
\begin{eqnarray*}
&&
\langle A (t, (R_1\bw)(t), \bw(t)) - \fb(t),
\bv - \bw(t) \rangle
+ \varphi (t, (R_3\bw)(t), \bw(t), \bv) \\[1mm]
&&\qquad \
- \, \varphi(t, (R_3\bw)(t), \bw(t), \bw(t))
+ \, j^0 (t, (R_4\bw)(t), M\bw(t); M\bv - M\bw(t)) \ge 0
\end{eqnarray*}
for all $\bv \in U$, a.e. $t \in (0,T)$.
\end{Problem}

We will apply Theorem~\ref{theorem1q} to prove the unique solvability of Problem~\ref{CONTACT2q}. To this end,
we will verify hypotheses $H(A)$, $H(f)$, $H(\varphi)_1$,
$H(j)_1$, $H(K)$, $H(M)$, $H(R)$, and $(H_0)$.
It follows from $H({\mathscr{A}})$(1) and (2)
that $H(A)$(a) holds.
By $H({\mathscr{A}})$(3), we have
\begin{eqnarray}
&&
|\langle A(t, \lambda, \bv), \bz \rangle | \le
(\sqrt{2} \, \| {\widetilde{a}}_0 (t) \|_{L^2(\Omega)}
+ {\widetilde{a}}_1 \| \bv \|) \, \| \bz \|
+ \| \lambda \|_E \| \bz \|
\end{eqnarray}
for all $t \in [0, T]$, $\lambda \in E$, $\bv$,
$\bz \in V$ which implies $H(A)$(b) with
$a_0(t) = \sqrt{2} \, \| {\widetilde{a}}_0 (t) \|_{L^2(\Omega)}$.
Based on $H({\mathscr{A}})$(2) and (3),
by a similar argument as used in~\cite[Theorem~7.3]{MOSBook},
we infer that $A(t, \lambda, \cdot)$ is continuous
for all $(t, \lambda) \in [0, T] \times E$ which ensures $H(A)$(c).
From $H({\mathscr{A}})$(4), we can assert that
$$
\langle A (t, \lambda, \bv_1)
- A(t, \lambda, \bv_2), \bv_1 - \bv_2 \rangle
\ge m_{\mathscr A} \| \bv_1 - \bv_2 \|^2
$$
for all $t \in [0, T]$, $\lambda \in E$,
$\bv_1$, $\bv_2 \in V$.
It is straightforward to obtain
$$
\| A (t, \lambda_1, \bv) - A(t, \lambda_2, \bv) \|_{V^*} \le \| \lambda_1 - \lambda_2\|_E
$$
for all $t \in [0, T]$, $\lambda_1$, $\lambda_2 \in E$,
$\bv \in V$.
From the last two inequalities, by Remark~\ref{REM2},
we deduce that $H(A)$(d) is satisfied with
$m_A = m_{\mathscr A}$.
Hence $H(A)$ is verified.
As already noted, under hypothesis $(H_1)$,
the functional $f$ defined by (\ref{fXXX}) satisfies
$\fb \in C(0, T; V^*)$. Hence $H(f)$ holds.

We will examine $H(\varphi)_1$.
By the convexity of the norm function,
it is immediate that
$\varphi(t, \eta, \bw, \cdot)$ is convex and lower semicontinuous for all $t \in [0, T]$,
$\eta \in Y$, and $\bw \in V$.
Let $t \in [0, T]$, $\eta_1 = (\eta_{11}, \eta_{12})$,
$\eta_2 = (\eta_{21}, \eta_{22}) \in Y$,
$\bw_1$, $\bw_2$, $\bv_1$, $\bv_2 \in V$.
We compute
\begin{eqnarray*}
&&\hspace{-0.5cm}
\varphi (t, \eta_1, \bw_1, \bv_2) -
\varphi (t, \eta_1, \bw_1, \bv_1) +
\varphi (t, \eta_2, \bw_2, \bv_1) -
\varphi (t, \eta_2, \bw_2, \bv_2) \\[1mm]
&&\hspace{-0.3cm}
= \int_{\Gamma_{3}}
(F_{b}(t, \eta_{11}) - F_{b}(t, \eta_{21}))
(\| \bv_{2\tau} \| - \| \bv_{1\tau} \|) \, d\Gamma
+ \int_{\Gamma_4}
(p(t, \eta_{12}) - p(t, \eta_{22}))(v_{2\nu}-v_{1\nu}) \, d\Gamma \\[1mm]
&&
+ \int_{\Gamma_{4}}
(\mu(t, \| \bw_{1\tau} \|) p(t, \eta_{12})
- \mu(t, \| \bw_{2\tau} \|) p(t, \eta_{22}))
\,
(\| \bv_{2\tau} \| - \| \bv_{1\tau} \|) \, d\Gamma
\\[1mm]
&&\ \
\le L_{F_b} \int_{\Gamma_{3}}
|\eta_{11} - \eta_{21}| \| \bv_{2\tau} - \bv_{1\tau} \|
\, d\Gamma
+ L_p \int_{\Gamma_4}
|\eta_{12} - \eta_{22}| |v_{2\nu}-v_{1\nu}| \, d\Gamma \\[1mm]
&&
+ \int_{\Gamma_{4}}
\Big(
\mu(t, \| \bw_{1\tau} \|) \, p(t, \eta_{12})
- \mu(t, \| \bw_{1\tau} \|) \, p(t, \eta_{22}))
+ \mu(t, \| \bw_{1\tau} \|) \, p(t, \eta_{22})\\[1mm]
&&\quad\qquad \qquad\qquad
- \, \mu(t, \| \bw_{2\tau} \|) \, p(t, \eta_{22}))
\Big) (\| \bv_{2\tau}- \bv_{1\tau} \|)\, d\Gamma
\\[1mm]
&&
\le L_{F_b} \| \gamma\| \| \eta_{11} - \eta_{21}\|_{L^2(\Gamma_3)} \| \bv_1 - \bv_2 \|
+ L_p \| \gamma\|
\|\eta_{12} - \eta_{22}\|_{L^2(\Gamma_{4})}
\| \bv_1 - \bv_2 \| \\[2mm]
&&
+ \, \mu_0 L_p \| \gamma\| \| \eta_{12} -\eta_{22}\|_{L^2(\Gamma_{4})} \| \bv_1 - \bv_2 \|
+
p^* L_\mu \| \gamma \|
\| \bw_{1} - \bw_{2\tau} \|_{L^2(\Gamma_{4};\real^d)}
\| \bv_1 - \bv_2 \| \, d\Gamma
\\[2mm]
&&
\le
(L_{F_b} + L_p + \mu_0 L_p) \| \gamma \|
\| \eta_1 - \eta_2 \|_Y \| \bv_1 - \bv_2\|
+
p^* L_\mu \| \gamma \|^2 \| \bw_1-\bw_2\| \| \bv_1 - \bv_2\|,
\end{eqnarray*}
which implies $H(\varphi)$(b) with
$\alpha_{\varphi} = p^* L_\mu \| \gamma \|^2$.
Subsequently,
from the inequality
\begin{eqnarray*}
&&\hspace{-0.5cm}
\varphi(t, \eta_1, \bw, \bv_1)
- \varphi(t, \eta_2, \bw, \bv_2)
\le \int_{\Gamma_{3}}
F_{b}(t, \eta_1) \| \bv_1 - \bv_2\| \, d\Gamma \\[1mm]
&&\hspace{-0.3cm}
+ \int_{\Gamma_4}
p(t, \eta_2) (v_{1\nu} -v_{2\nu})\, d\Gamma
+
\int_{\Gamma_4}
\mu(t, \| \bw_\tau \|) \,
p(t, \eta_2)  \| \bv_1 -\bv_{2} \| \, d\Gamma
\\[1mm]
&&
\le
c \, (1 + \| \eta\|_Y) \| \gamma \| \|\bv_1 -\bv_2\|
+ p^* \int_{\Gamma_{4}} \| \bv_1-\bv_2\| \, d\Gamma
+ \mu^* p^*
\int_{\Gamma_{4}} \| \bv_1-\bv_2\| \, d\Gamma \\[1mm]
&&\quad
\le c \, (1 + \| \eta \|_Y)\, \| \bv_1-\bv_2\|
\end{eqnarray*}
for all $t \in [0, T]$, $\eta = (\eta_{1}, \eta_{2})$,
$\bw$, $\bv_1$, $\bv_2 \in V$ with a constant $c > 0$,
we deduce $H(\varphi)$(c).

Let $\bv \in V$, $t_n \to t$ in $[0, T]$,
$\eta_n = (\eta_{1n}, \eta_{2n}) \to
(\eta_{1}, \eta_{2}) = \eta$ in $Y$ and
$\bw_n \rightharpoonup \bw$ in $V$.
By the compactness of the trace operator, we may suppose
that $\bw_n \to \bw$ in $L^2(\Gamma;\real^d)$.
We use $H(F_b)$(2),\,(3),
$H(p)$(2),\,(5),
$H(\mu)$(2),\,(4) and from the Lebesgue-dominated convergence theorem, see, e.g.,~\cite[Theorem~2.38]{MOSBook},
we obtain
\begin{eqnarray*}
&&\hspace{-0.5cm}
\varphi (t_n, \eta_n, \bw_n, \bv) -
\varphi (t_n, \eta_n, \bw_n, \bw_n)
- \varphi (t, \eta, \bw, \bv)
+ \varphi (t, \eta, \bw, \bw) \\[1mm]
&&
\le
\int_{\Gamma_{3}}
|F_b(t_n, \eta_{1n}) - F_b(t, \eta_1)|
\, \| \bv_{\tau} \| \, d\Gamma
+
\int_{\Gamma_{4}}
|p(t_n, \eta_{2n}) - p(t, \eta_2)| | v_{\nu} | \, d\Gamma
\\[1mm]
&&
+
\int_{\Gamma_{4}}
| \mu(t_n, \| \bw_{n\tau}\|) \, p(t_n, \eta_{2n})
- \mu(t, \| \bw_\tau \|) \, p(t, \eta_2)|
\, \| \bv_{\tau} \| \, d\Gamma \\[1mm]
&&
\le
\int_{\Gamma_{3}}
|F_b(t, \eta_{1}) \| \bw_\tau \|
- F_b(t_n, \eta_{1n}) \| \bw_{n\tau} \| | \, d\Gamma
+
\int_{\Gamma_{4}}
|p(t, \eta_{2}) w_\nu - p(t_n, \eta_{2n}) w_{n\nu} | \, d\Gamma \\[1mm]
&&
+
\int_{\Gamma_{4}}
| \mu(t, \| \bw_{\tau}\|) \, p(t, \eta_{2})\,\| \bw_\tau \|
- \mu(t_n, \| \bw_{n\tau} \|) \, p(t_n, \eta_{2n}) \,
\| \bw_{n\tau} \| | \, d\Gamma \to 0.
\end{eqnarray*}
Hence, we deduce condition $H(\varphi)_1$(d).
This completes the verification of $H(\varphi)_1$.

We next claim that condition $H(j)_1$ holds.
Analogously as in~\cite[Theorem~13]{hms},
we can prove that condition $H(j)$ is satisfied with
$m_j = k^* m_{j_\nu}\| \gamma \|^2$ and
$m_1 = c_0 L_k \| \gamma \|$.
For the proof of (\ref{EE1}), let
$\bv \in V$,
$t_n \to t$ in $[0, T]$,
$\zeta_n \to \zeta$ in $Z$ and
$\bw_n \rightharpoonup \bw$ in $V$.
By virtue of~\cite[Proposition~3.23(ii)]{MOSBook} and
$H(j_\nu)$, we have
\begin{equation}\label{GG44}\limsup j_\nu^0 (v_\nu; v_\nu-w_{n\nu}) \le
j_\nu^0(v_\nu; v_\nu-w_\nu).
\end{equation}
Using~\cite[Propositions~3.23(iii) and 3.47(iv)]{MOSBook}, it follows
\begin{eqnarray}
&&\label{GG55}
\hspace{-0.5cm}
j^0(t_n, \zeta_n, v_\nu; v_\nu-w_{n\nu}) \le
\int_{\Gamma_{3}} k(t_n, \zeta_n)
j_\nu^0 (v_\nu; v_\nu-w_{n\nu}) \, d\Gamma
\\[1mm]
&&\nonumber
\le \int_{\Gamma_{3}} |k(t_n, \zeta_n) - k(t, \zeta)| \,
|\partial j_\nu (v_\nu)| \, | v_\nu - w_{n\nu}| \, d\Gamma
+ \int_{\Gamma_{3}}
k(t, \zeta) j_\nu^0 (v_\nu; v_\nu - w_{n\nu})
\, d\Gamma.
\end{eqnarray}
We pass to the upper limit in (\ref{GG55}),
combine with (\ref{GG44}), and use the Fatou lemma
and $H(k)$ to get
$$
\limsup j^0(t_n, \zeta_n, Mv; Mv-Mw_n)
\le \int_{\Gamma_{3}}
k(t, \zeta) j_\nu^0 (Mv; Mv - Mw) \, d\Gamma.
$$
Finally, in the latter, we exploit the regularity hypothesis $H(j_\nu)$(1), and deduce $H(j)_1$.

It is clear that the set $K= U$ is a closed and convex subset of $V$ with $\bzero \in K$, i.e., $H(K)$ holds.
By the properties of the (normal) trace operator, it is
clear that $H(M)$ holds.
Furthermore,
by hypotheses $H({\mathscr B})$ and $H({\mathscr C})$,
we obtain that
the operators $I$, $R_1$, $R_3$ and $R_4$ satisfy $H(R)$,
for details see~\cite[Theorem~14.2]{SMO}.
The smallness condition $(H_0)$ is a consequence of
(\ref{CONTACT_smallness}).

We have verified all hypotheses of Theorem~\ref{theorem1q}.
Then, we deduce that Problem~\ref{CONTACT2q}
has a unique solution $\bw \in L^2(0, T; U)$.
We define a function $\bu \colon (0, T) \to V$ by
$$
\bu(t) = \bu_0 + \int_0^t \bw (s)\,ds
\ \ \mbox{for} \ \ t \in [0, T].
$$
By the regularity hypothesis $H(j_\nu)$(1) it follows that
$\bw \in L^2(0, T; V)$ solves Problem~\ref{CONTACT2q} if and only if $\bu \in C(0, T; V)$ is the unique solution to Problem~$\ref{CONTACT1q}$.
In conclusion, Problem~$\ref{CONTACT1q}$ has a unique solution
$\bu \in C(0, T; V)$ with $\bu' \in L^2(0, T; V)$ and
$\bu'(t) \in U$ for a.e. $t \in (0, T)$.
This completes the proof.
\hfill$\Box$


\begin{thebibliography}{99}


\bibitem{CLM1}
S. Carl, V.K. Le, D. Motreanu,
{\it Nonsmooth Variational Problems and their Inequalities. Comparison Principles and Applications}, Springer Monogr. Math.,
Springer, New York, 2007.






\bibitem{cla}
F.H. Clarke, {\it Optimization and Nonsmooth Analysis}, Wiley, New York, 1983.


\bibitem{DMP1} Z. Denkowski, S. Mig\'{o}rski, N.S. Papageorgiou, {\it An Introduction to Non\-li\-ne\-ar Analysis: Theory}, Kluwer Academic/Plenum Publishers, Boston, Dordrecht, London, New York, 2003.









\bibitem{hms}
W. Han, S. Mig\'{o}rski, M. Sofonea, Analysis of a general dynamic history--dependent
variational--hemivariational inequality, {\it Nonlinear Analysis: Real World Applications} {\bf 36} (2017), 69--88.

\bibitem{HS}
W. Han, M. Sofonea, {\it Quasistatic  Contact Problems in Viscoelasticity and Viscoplasticity},
Studies in Advanced Mathematics {\bf 30}, Americal Mathematical Society, Providence, RI--International Press, Somerville, MA, 2002.


\bibitem{KS}
D. Kinderlehrer, G. Stampacchia, {\em An Introduction to Variational Inequalities and their Applications}, Classics in Applied Mathematics {\bf 31}, SIAM, Philadelphia, 2000.

\bibitem{KULIG}
A. Kulig, S. Mig\'orski,
Solvability and continuous dependence results for second order nonlinear inclusion with Volterra-type operator, {\it Nonlinear Analysis: Theory, Methods \& Applications} {\bf 75} (2012), 4729--4746.



\bibitem{MigJOTA}
S. Mig\'orski, Optimal control of history--dependent evolution inclusions with applications to frictional contact, {\it Journal of Optimization Theory and Applications} {\bf 185} (2020), 574--596.

\bibitem{Mig2021}
S. Mig\'orski,
A class of history--dependent systems of evolution inclusions with applications, {\it Nonlinear Analysis: Real World Applications} {\bf 59} (2021), 103246.

\bibitem{MigBai}
S. Mig\'orski, Y.R. Bai,
Well-posedness of history--dependent evolution
inclusions with applications, {\it Zeitschrift f\"ur angewandte Mathematik und Physik ZAMP} {\bf 70} (2019), 114.




\bibitem{MOS13}
S. Mig\'orski, A. Ochal, M. Sofonea, History--dependent
subdifferential inclusions and hemivariational inequalities in contact mechanics, {\it Nonlinear Analysis: Real World Applications} {\bf 12} (2011), 3384--3396.

\bibitem{MOSBook}
S. Mig\'orski, A. Ochal, M. Sofonea,
{\it Nonlinear Inclusions and Hemivariational Inequalities. Models and Analysis of Contact Problems},
Advances in Mechanics and Mathematics \textbf{26},
Springer, New York, 2013.

\bibitem{MOS18}
S. Mig\'orski, A. Ochal, M. Sofonea, History--dependent
variational--hemi\-va\-ria\-tio\-nal inequalities in contact mechanics, {\it Nonlinear Analysis: Real World Applications} {\bf 22} (2015), 604--618.

\bibitem{MOS30}
S. Mig\'orski, A. Ochal, M. Sofonea,
A class of variational-hemivariational inequalities in reflexive Banach spaces, {\it Journal of Elasticity} {\bf 127} (2017), 151--178.


\bibitem{MSZ2019}
S. Mig\'orski, M. Sofonea, S.D. Zeng, Well-posedness of history--dependent sweeping processes,
{\it SIAM Journal on Mathematical Analysis} {\bf 51}  (2019), 1082--1107.

\bibitem{MBZ}
S. Mig\'orski, B. Zeng, A new class of history-dependent evolutionary variational-hemivariational inequalities with unilateral constraints,
{\it Applied Mathematics and Optimization} {\bf 84} (2021), 2671--2697.

\bibitem{MigZENG}
S. Mig\'orski, S.D. Zeng, Rothe method and numerical analysis for history--dependent hemivariational inequalities with applications to contact mechanics, {\it Numerical Algorithms} {\bf 82} (2019), 423--450.


\bibitem{MOTSOF}
D. Motreanu, M. Sofonea,
Evolutionary variational inequalities arising in quasistatic frictional contact problems for elastic materials, {\it Abstract and Applied Analysis} {\bf 4} (1999), 255--279.

\bibitem{NP}
Z. Naniewicz, P.D. Panagiotopoulos,
{\it Mathematical Theory of Hemivariational Inequalities
and Applications}, Marcel Dekker, Inc., New York, Basel, Hong Kong, 1995.

\bibitem{Oden}
J.T. Oden, J.A.C. Martins,
Models and computational methods for dynamic friction phenomena, {\it Computer Methods in Applied Mechanics and Engineering}
{\bf 52} (1985), 527--634.



\bibitem{P2}
P.D. Panagiotopoulos, {\it Inequality Problems in Mechanics and Applications. Convex and Nonconvex
Energy Functions}, Birkh\"auser, Basel, 1985.

\bibitem{P}
P.D. Panagiotopoulos, {\it Hemivariational Inequalities,
Applications in Mechanics and Engineering}, Springer-Verlag, Berlin, 1993.

\bibitem{SST}
M. Shillor, M. Sofonea, J.J. Telega, {\it Models and Analysis of Quasistatic Contact}, Lect. Notes Phys. {\bf 655}, Springer, Berlin, Heidelberg, 2004.

\bibitem{Sofonea}
M. Sofonea,
Optimal control of a class of variational-hemivariational inequalities in reflexive Banach spaces, {\it Applied Mathematics \& Optimization} {\bf 79} (2019), 621--646.

\bibitem{SHM2015}
M. Sofonea, W. Han, S. Mig\'orski, Numerical analysis of history--dependent va\-ria\-tio\-nal-hemivariational
inequalities with applications to contact problems,
{\it European Journal of Applied Mathematics} {\bf 26} (2015), 427--452.

\bibitem{SMatei2011}
M. Sofonea, A. Matei,
History--dependent quasi-variational inequalities arising in
contact mechanics, {\it European Journal of Applied Mathematics} {\bf 22} (2011), 471--491.

\bibitem{SMatei}
M. Sofonea, A. Matei, {\em Mathematical Models in Contact Mechanics}, London Mathematical Society Lecture Notes Series {\bf 398}, Cambridge University Press, Cambridge, 2012.

\bibitem{SM2016}
M. Sofonea, S. Mig\'orski, A class of history-dependent variational-hemivariational inequalities,
{\it Nonlinear Differential Equations and Applications}
{\bf 23} (2016), 38.

\bibitem{SM2}
M. Sofonea, S. Mig\'orski,
{\it Variational--Hemivariational Inequalities with Applications},
Chapman \& Hall/CRC, Monographs and Research Notes in Mathematics, Boca Raton, 2018.

\bibitem{SMO}
M. Sofonea, S. Mig\'orski, A. Ochal,
Two history--dependent contact problems, in: W. Han, S. Mig\'orski, M. Sofonea (Eds.),
{\it Advances in Variational and Hemivariational Inequalities: Theory, Numerical Analysis, and Applications}, in: Advances in
Mechanics and Mathematics Series, vol. 33,
Springer, 2015, 355--380, Chapter 14.

\bibitem{SPatrulescu}
M. Sofonea, F. P\v{a}trulescu,
Penalization of history--dependent variational inequalities, {\it European Journal of Applied Mathematics} {\bf 25} (2014), 155--176.

\bibitem{SXiao2015}
M. Sofonea, Y.B. Xiao,
Fully history--dependent quasivariational inequalities in contact mechanics, {\it Applicable Analysis} {\bf 95} (2016),  2464--2484.

\bibitem{Xiao2019}
Y.B. Xiao, M. Sofonea,
Generalized penalty method for elliptic
variational-he\-mi\-va\-ria\-tional inequalities,
{\it Applied Mathematics \& Optimization} {\bf 83} (2021),  789--812.


\end{thebibliography}
\end{document}